\documentclass[reqno, 12pt]{amsart}

\usepackage{amsmath,amssymb,amsthm}
\usepackage[T1]{fontenc}
\usepackage{lmodern}
\usepackage[linktocpage,colorlinks,linkcolor=magenta,citecolor=magenta,urlcolor=black,pagebackref=false]{hyperref}
\usepackage[noEnd=false]{algpseudocodex}
\usepackage{enumitem}
\usepackage[all]{xy}
\usepackage{microtype}
\usepackage{xcolor}
\usepackage{cleveref}  

\newtheorem{thm}{Theorem}[section]
\newtheorem{lem}[thm]{Lemma}
\newtheorem{cor}[thm]{Corollary}
\newtheorem{prop}[thm]{Proposition}

\theoremstyle{definition}
\newtheorem{defn}[thm]{Definition}

\newtheorem{rem}[thm]{Remark}
\newtheorem{ex}[thm]{Example}
\newtheorem{algo}[thm]{Algorithm}
\newtheorem{proc}[thm]{Procedure}

\let\oldphi\phi
\let\phi\varphi
\let\varphi\oldphi
\let\oldepsilon\epsilon
\let\epsilon\varepsilon
\let\varepsilon\oldepsilon
\newcommand{\eps}{\epsilon}

\newcommand{\mbb}[1]{\ensuremath{\mathbb{#1}}}
\newcommand{\Z}{\mbb{Z}}
\newcommand{\N}{\mbb{N}}

\newcommand{\FF}{\mbb{F}}

\newcommand{\mbf}[1]{\ensuremath{\mathbf{#1}}}

\newcommand{\Pb}{\mbf{P}}
\newcommand{\Fb}{\mbf{F}}
\newcommand{\Gb}{\mbf{G}}
\newcommand{\Ib}{\mbf{I}}
\newcommand{\Jb}{\mbf{J}}
\newcommand{\Mb}{\mbf{M}}
\newcommand{\Nb}{\mbf{N}}
\newcommand{\Ab}{\mbf{A}}

\newcommand{\fm}{\mathfrak{m}}

\newcommand{\mcal}[1]{\ensuremath{\mathcal{#1}}}
\newcommand{\Ac}{\mathcal{A}}
\newcommand{\Bc}{\mathcal{B}}

\newcommand{\Fc}{\mcal{F}}
\newcommand{\Gc}{\mcal{G}}
\newcommand{\Ic}{\mcal{I}}
\newcommand{\Jc}{\mcal{J}}
\newcommand{\Mcc}{\mathcal{M}}
\newcommand{\Nc}{\mathcal{N}}
\newcommand{\Pc}{\mathcal{P}}
\newcommand{\Sc}{\mathcal{S}}

\DeclareMathOperator{\Inc}{Inc}
\DeclareMathOperator{\Incn}{Inc (\N)}
\DeclareMathOperator{\Sym}{Sym}
\DeclareMathOperator{\Hom}{Hom} 
\DeclareMathOperator{\Mod}{Mod}
\DeclareMathOperator{\rank}{rank}
\DeclareMathOperator{\Tor}{Tor}

\newcommand{\tn}[1]{\textnormal{#1}}
\let\hom\relax\newcommand{\hom}{\tn{Hom}}
\newcommand{\im}{\tn{im}}

\newcommand{\id}{\tn{id}}

\newcommand{\la}{\langle}
\newcommand{\ra}{\rangle}
\newcommand{\FI}{\tn{FI}}
\newcommand{\OI}{\tn{OI}}

\renewcommand{\mod}{\tn{-Mod}}

\renewcommand{\l}{\ell}

\numberwithin{equation}{section}
\algrenewcommand\algorithmicrequire{\textbf{Input:}}
\algrenewcommand\algorithmicensure{\textbf{Output:}}

\title{Equivariant Free Resolutions of  Sequences of Symmetric Modules}

\author{Michael Morrow}
\address{Department of Mathematics, University of Kentucky, 715 Patterson Office Tower, Lexington,
KY 40506 USA}
\email{michael.morrow@uky.edu}

\author{Uwe Nagel}
\address{Department of Mathematics, University of Kentucky, 715 Patterson Office Tower, Lexington,
KY 40506 USA}
\email{uwe.nagel@uky.edu}

\thanks{The second author was partially supported by Simons Foundation grant \#636513.  \\
}

\begin{document}

\begin{abstract} Given a sequence of related modules $M_n$ over a sequence of related Noetherian polynomial rings, where each $M_n$ is a representation of the symmetric group on $n$ letters, one may ask how to simultaneously compute an equivariant free resolution of each $M_n$. In this article, we address this question. Working in the setting of FI-modules over a Noetherian polynomial FI-algebra, we provide an algorithm for computing syzygies and FI-equivariant differentials. As an application, we show how this result can be used to compute truncations of equivariant free resolutions of ideals in polynomial rings in infinitely many variables that are invariant under actions of the monoid of strictly increasing maps or of permutations. The free modules occurring in such a free resolution are finitely generated up to symmetry. 
\end{abstract}

\maketitle

\tableofcontents


\section{Introduction}
\label{section:intro} 

The development of a theory of FI-modules over an FI-algebra has been motivated by questions in algebraic statistics (see, e.g., \cite{BD, Draisma-factor, HS}) and the success of methods used to establish representation stability (see, e.g., \cite{CE, CEF, CEFN, SS-14}). 
Intuitively, an FI-module  may be considered as a sequence of related symmetric $\Ab_n$-modules $\Mb_n$, where the rings $\Ab_n$ and modules $\Mb_n$ depend on an integer $n \ge 0$. If $M$ is a finitely generated  $P$-module, where $P$ is a polynomial ring over a field $K$ in finitely many variables, then $M$ admits a finite free resolution over $P$ that can be computed algorithmically (see, e.g., \cite{Eis95}). One goal of this article is to establish an analogous result for a sequence of related symmetric modules $\Mb_n$ over polynomial rings $\Pb_n$. It provides simultaneously a (truncated) free resolution of  each $\Mb_n$ over $\Pb_n$ whose differentials are equivariant with respect to an action of the symmetric group $\Sym(n)$ on $n$ letters. 

More precisely, fix some integer $c \ge 0$ and denote by $\Pb_n$ a polynomial ring $K[x_{i, j} \; \mid \; i \in [c], j \in [n]]$  in $c n$ variables over a field $K$, where $[n]$ is the interval $\{1,\ldots,n\}$.  The rings $\Pb_n$ can be used to define a functor $\Pb$ from the category FI to the category of $K$-algebras with $\Pb_n = \Pb([n])$, where FI is the category of finite sets and injections.  An FI-module over $\Pb$ is a functor $\Mb$ from FI to the category of $K$-modules with suitable compatibility properties (see \Cref{def:FI-module} for details). In particular, each $\Mb_n = \Mb([n])$ is a $\Pb_n$-module and a $K$-linear representation of $\Sym (n)$. In \cite{NR19}, it is shown that any finitely generated FI-module $\Mb$ over $\Pb$ is Noetherian, and so $\Mb$ admits an exact sequence of finitely generated $\Pb$-modules 
\[
\Fb^{\bullet}: \hspace*{2cm} \cdots \to \Fb^2 \to \Fb^1 \to \Fb^0 \to \Mb \to 0, 
\]
where each $\Fb^i$ is a free FI-module over $\Pb$. We call such a sequence a \emph{free $\FI$-resolution} of $\Mb$. For every integer $n \ge 0$, it immediately gives a free resolution 
\[
\Fb^{\bullet}_n: \hspace*{2cm} \cdots \to \Fb^2_n \to \Fb^1_n \to \Fb^0_n \to \Mb_n \to 0
\]
of $\Mb_n$ over the polynomial ring $\Pb_n$ that is equivariant with respect to the action of $\Sym (n)$ (see \Cref{prop:FI gives symmetry}). In this article, we present a finite algorithm that computes the first $s$ steps of a free FI-resolution of $\Mb$ for any given integer $s \ge 0$ (see \Cref{proc:freefires}). It utilizes an analogous result for OI-modules in \cite{MN} that relies on the  use of  OI-Gr\"obner bases. Note that the definition of a Gr\"obner basis requires a suitable order, which is in conflict with the symmetry induced by the FI-structure. Thus, the methods in \cite{MN} do not directly apply to FI-modules. 

If $\Mb$ is a graded FI-module over $\Pb$ endowed with the standard grading where every variable has degree one,  then $\Mb$ admits a graded free FI-resolution where each of the free modules $\Fb^i$ has minimal rank, and such a minimal resolution is unique. This follows from work in \cite{FN21}. A variant of our algorithm determines truncations of such a minimal resolution $\Fb^{\bullet}$ of $\Mb$. The corresponding graded free resolution $\Fb^{\bullet}_n$ of $\Mb_n$ is not necessarily minimal as a resolution over $\Pb_n$ in the classical sense. However, it can be algorithmically modified to obtain an equivariant  minimal free resolution of $\Mb_n$ (see \Cref{prop:prune res}). 

Recently, symmetric ideals in a polynomial ring $K[X_\N] = K[x_1, x_2,\ldots]$ have received considerable attention (see, e.g., \cite{NS-20, NS-21}).   Notice that, for any $c \ge 0$,  a polynomial ring $\Pc = K[X_{c \times \N} = K[x_{i, j} \; \mid \; i \in [c], j \in \N]$  is isomorphic to the colimit of the rings $\Pb_n$ and that $\Sym = \bigcup_{n \in \N} \Sym (n)$ acts on $\Pc$ by permuting column indices of the variables. In order to leverage this connection, we introduce more generally Sym-algebras and Sym-modules over them. Our prime example of a Sym-algebra is $\Pc$. The colimit of a finitely generated FI-module $\Mb$ over $\Pb$ is a $\Pc$-module $\Mcc$ that is finitely generated up to symmetry. Using functoriality of colimits, we show that, for any fixed integer $s \ge 0$, there is a finite algorithm to compute  the first $s$ steps of a free resolution of any symmetric ideal $\Ic$ of $\Pc$, 
\[
 \Fc_s \to \Fc_{s-1} \to \cdots \to  \Fc_1 \to \Ic \to 0 
\]
with Sym-equivariant differentials, where each $\Fc_i$ is a  free Sym-module over $\Pc$ that is finitely generated up to symmetry. 

This article is organized as follows: In the next section, we review needed concepts for FI- and OI-modules, where OI denotes the category of totally ordered finite sets and order-preserving injections. We also include a discussion of minimal and equivariant free resolutions over polynomial rings in finitely many variables. In \Cref{section:FI-res}, we present the first algorithm to compute a truncated free FI-resolution of a finitely generated FI-module over a Noetherian polynomial FI-algebra. 

The mentioned polynomial rings $\Pb_n$ can also be used to define an OI-algebra denoted $\Pb^\OI$. Its colimit admits an action of the monoid $\Inc$ of strictly increasing maps $\N \to \N$. Building on initial work in \cite{HS, NR19}, we develop the foundations of a theory of $\Inc$-modules over an $\Inc$-algebra in \Cref{section:Inc-resolutions}. We show that any submodule of a finitely generated free $\Inc$-module over $\Pb^\OI$ admits an $\Inc$-equivariant free resolution (see \Cref{thm:Inc resolutions}). Furthermore, combined with results in \cite{MN}, we prove that the first $s$ steps of such a resolution can be computed algorithmically (see \Cref{proc:Inc-free res}). We also present several explicit examples. 

In \Cref{section:Sym-resolutions}, we discuss symmetric modules over $\Pc$, considered as a $\Sym$-algebra. To this end, we introduce more generally $\Sym$-modules over a $\Sym$-algebra. The results are analogous to those in \Cref{section:Inc-resolutions}. In \Cref{thm:Sym resolutions},  we establish, for any submodule of a finitely generated free $\Sym$-module over $\Pc$, the existence of a free resolution with $\Sym$-equivariant differentials. We show that truncations of such a resolution can be computed in finitely many steps (see \Cref{proc:Sym-free res}).


\section{Preliminaries}
\label{section:prelim}

We discuss background on OI- and FI-modules that is used later on. Though we begin with some considerations  on minimal and equivariant free resolutions in the classical setting of Hilbert's Syzygy Theorem. 


Denote by  $P$ a graded polynomial ring in finitely many variables over a field $K$ with a standard grading, i.e., every variable has degree one. Let $M = \oplus_{j \in \Z} [M]_j$ be a finitely generated graded module over $P$. Then it is well-known that $M$ admits a finite graded minimal free resolution 
\[
0 \to F_s \stackrel{\phi_s}{\longrightarrow} F_{s-1} \to \cdots \to F_1 \stackrel{\phi_1}{\longrightarrow} F_0 \to M \to 0  
\]
with finitely generated free $P$-modules $F_i$ and minimal degree-preserving  homomorphisms~$\phi_i$. Recall that $\phi_i$ is said to be \emph{minimal} if its image is contained in $\fm F_{i-1}$, where $\fm$ is the maximal ideal of $P$ generated by all its variables. A graded minimal free resolution is unique up to isomorphisms of complexes. 

Assume the symmetric group on $n$ letters $\Sc = \Sym (n)$ acts $K$-linearly on $P$ and $M$ such that it preserves degrees, induces ring homomorphisms on $P$ and is \emph{compatible} with the $P$-module structure, i.e., $\sigma (pm) = (\sigma \cdot p) (\sigma \cdot m)$ for any $\sigma \in \Sc, p \in P$ and $m \in M$. Then $M$ admits an \emph{equivariant} graded minimal free resolution, i.e., the maps $\phi_i$ satisfy additionally $\phi_i (\sigma \cdot f) = \sigma \cdot (\phi_i (f)$. It seems part of the folklore that one can adjust any given equivariant graded free resolution of $M$ to obtain an equivariant graded minimal free resolution (in the non-modular case). For lack of a suitable reference, we include a proof of this result here. 

\begin{prop}
     \label{prop:prune res}
If $M$ has a graded free resolution $\FF$ over $P$, 
\[
0 \to F_s \stackrel{\phi_s}{\longrightarrow} F_{s-1} \to \cdots \to F_1 \stackrel{\phi_1}{\longrightarrow} F_0  \stackrel{\phi_1}{\longrightarrow}  M \to 0,    
\]
then one has:
\begin{itemize}

\item[{\rm (a)}] One can modify $\FF$ to produce a minimal graded free resolution of $M$ in finitely many steps. 

\item[{\rm (b)}] If $\FF$ is equivariant and the characteristic of $K$ is zero or does not divide $n$ then one can modify $\FF$ to obtain an equivariant  minimal graded free resolution of $M$ in finitely many steps. 

\end{itemize}

\end{prop}

\begin{proof}
We first consider (a). It is well-know that a graded free resolution over $P$ is minimal if and only if one cannot split off a trivial complex. A complex is called \emph{trivial} if it is isomorphic to a complex of the form 
\[
0 \to F \stackrel{\id}{\longrightarrow} F \to 0, 
\]
where $F$ is a finitely generated graded free $P$-module. Thus, it is enough to establish the following claim and to apply it repeatedly. 
\smallskip

\emph{Claim:} Consider a complex of  finitely generated graded $P$-modules 
\[
Q \stackrel{\gamma}{\longrightarrow} F \stackrel{\phi}{\longrightarrow} G \stackrel{\beta}{\longrightarrow} N  
\]
with free $P$-modules $F, G, Q$, which is exact at $F$ and $G$. If $\phi$ is not a minimal homomorphism then one can split off a trivial complex to obtain an analogous complex 
\[
Q \stackrel{\gamma'}{\longrightarrow} F' \stackrel{\phi'}{\longrightarrow} G' \stackrel{\beta'}{\longrightarrow} N,   
\]
where $\rank F' = \rank F - 1$ and $\rank G' = \rank G - 1$. 
\smallskip 

Indeed, choose (ordered) bases $\{f_1,\ldots,f_s\}$ and $\{g_1,\ldots,g_t\}$ of $F$ and $G$, respectively. There are homogeneous $p_{i, j} \in P$ such that 
\[
\phi (f_i) = \sum_{j=1}^t p_{i, j} g_j. 
\]
Since $\phi$ is not minimal one of the coefficients $p_{i, j}$ must be unit. Switching bases vectors, we may assume $a = p_{1,1} \in K$ is a unit. Thus, the set $\{g_1',\ldots,g_t'\}$ defined by   
\[
    g_{1}' = \phi(f_{1}) \quad \text{ and } \quad g_{j}' = g_{j} \;
    \text{ if } 2 \le j \le t
\]
is  another basis of $G$ . Similarly, we get a new basis  $\{f_1',\ldots,f_s'\}$ of $F$ by setting
\[
    f_{1}' = f_{1} \quad \text{ and }  \quad  f_{i}' = f_{i} - a^{-1} \cdot p_{i, 1} f_1 \; \text{ if } 2 \le i \le s.  
\]
Denote by $F'$ and $G'$ the free $P$-submodules of $F$ and $G$ generated by $\{f_2',\ldots,f_s'\}$ and $\{g_2',\ldots,g_t'\}$, respectively. Notice that $\phi (F') \subseteq G'$ because
\[
\phi (f'_i) = \sum_{j=2}^t [ p_{i, j} - a^{-1} p_{1, j}  p_{i, 1} ] g_j \in G'. 
\] 
Hence, we can define a $P$-module homomorphism  $\phi' \colon F' \to G'$ by setting $\phi' (f) = \phi (f)$ for any $f \in F'$. Furthermore, we define a homomorphism of free rank one $P$-modules  $\alpha \colon \langle f_1\rangle \to \langle g_1' \rangle$ by putting $\alpha (f_1) = \phi (f_1) = g_1'$. It follows that $\phi = \alpha \oplus \phi'$. Since $\alpha$ is an isomorphism we can split off the trivial complex 
\[
\langle f_1\rangle \stackrel{\alpha'}{\longrightarrow} \langle g_1' \rangle
\]
and obtain a complex
\[
Q \stackrel{\gamma'}{\longrightarrow} F' \stackrel{\phi'}{\longrightarrow} G' \stackrel{\beta'}{\longrightarrow} N,   
\]
where $\beta'$ is the restriction of $\beta$ to $G'$ and $\gamma'$ is defined by $\gamma' (q)  = \gamma (q)$ for any $q \in Q$. This new complex is exact at $F'$ and $G'$, which completes the argument for establishing the Claim and (a). 

Second, we consider Assertion (b). By the assumptions on the characteristic of $K$, the symmetric group $\Sc$ is linearly reductive. By part (a), we may assume that  we used the given equivariant free resolution to get a minimal free resolution. It remains to show how to adjust the action of $\Sc$ on the free modules $F_i$ to obtain 
an equivariant minimal free resolution. Indeed, this can be achieved by applying the following method to the homomorphisms $F_i \to \phi (F_i)$ induced by $\phi_i$ for $i = 0, 1,\ldots,s$ in this order. 

Consider a graded  surjective $P$-module homomorphism $\phi \colon F \to N$ of finitely generated graded $P$-modules, where $F$ is a free module and $\phi$ is a minimal homomorphism, that is, the images $\phi (f_i)$ of the elements of a basis $\{f_1,\ldots,f_s\}$ of $F$ 
form a minimal generating set  of $N$. We will define a new group action of $\Sc$ on $F$ that is again degree-preserving and compatible with the $P$-module structure and a minimal surjective $P$-module homomorphism $\phi' \colon F \to N$ such that $\phi'$ is equivariant with respect to the new group action and satisfies $\ker \phi' = \ker \phi$.

Since $N$ is $\Sc$-invariant so is the module $N/\fm N$. Moreover, it decomposes into a direct sum 
$N/\fm N = \oplus_{j \in \Z} [N/\fm N]_j$, in which only finitely many of the degree components $[N/\fm N]_j$ are nonzero. These degree components are $\Sc$-invariant as well since the given group action preserves degrees. For any degree $j$ such that $[N/\fm N]_j \neq 0$, choose a set $G_j = \{n_{j, 1},\ldots,n_{j, k_j}\}$ of elements of $[N]_j$ such that their images form a basis of the $K$-vector space $[N/\fm N]_j$ and,  
for each $\sigma \in \Sc$ and any $n_{j, k} \in G_j$, the element $\sigma \cdot n_{j, k}$ is a linear combination of the elements in $G_j$, that is, 
\begin{align}
    \label{eq:action on N}
\sigma \cdot n_{j, k} = \sum_{l = 1}^{k_j} \lambda_{l} n_{j, l} \quad \text{ with }  \lambda_{l} \in K  \text{ depending on $\sigma, j$ and $k$}
\end{align}
(see, e.g., Theorem 8 in Section 2.6 of \cite{Serre}).  

By Nakayama's Lemma, the disjoint union $\bigcup_j G_j$ is a minimal generating set of $N$. Thus its cardinality is $s = \rank F$. Hence, we can reindex the basis elements of $F$ and define a homomorphism 
\[
\phi' \colon F \to N \quad \text{ by setting } \; \phi' (f_{j, k}) = n_{j, k}. 
\]
It is surjective and minimal. 

Using that both sets $\{\phi (f_{j, k}) \; | \; j, k\}$ and $\{ n_{j, k}  \; | \; j, k \}$ are minimal generating sets of $N$, 
 there is an isomorphism $\alpha \colon N \to N$ with $\alpha (\phi (f_{j, k})) = n_{j, k}$. It follows that $\ker \phi' = \ker \phi$. 

Since the coefficients $\lambda_{l} \in K$ in Equation~\eqref{eq:action on N}  are uniquely determined given $\sigma, j$ and $k$  we may define a new action of $\Sc$ on $F$ by setting 
\[
\sigma* (p f_{j, k}) = (\sigma \cdot p)  \sum_{l = 1}^{k_j} \lambda_{l} f_{j, l}  
\]
if $\sigma \cdot n_{j, k} = \sum_{l = 1}^{k_j} \lambda_{l} n_{j, l} $ and $p \in S$ and extending this additively. This does indeed give a group action on $F$. 
By definition, the map $\phi$ is equivariant with respect to this action of $\Sc$, which completes the argument. 
\end{proof}

\begin{rem}
   \label{rem:irr decomposition} 
Using an $\Sc$-equivariant graded minimal free resolution of $M$, one can algorithmically compute  decompositions of  graded components of $M$ and its Tor modules $\Tor^i_P (M, K)$ as direct sum of irreducible representations of $\Sc$ algorithmically (see \cite{Galetto-22, Galetto-23}). 
\end{rem}


Next, we recall needed background on OI- and FI-modules and fix notation. Let $K$ be a commutative ring with identity. 

\begin{defn}[{\cite[Definition 2.4]{NR19}}]
   \label{def:OI, FI algebra}
(i) Denote by FI  the category of finite sets and  injective maps. An \emph{\FI-algebra} over $K$ is a (covariant) functor $\Ab$ from FI to the category of commutative, associate, unital $K$-algebras. 

(ii) The subcategory of totally ordered finite sets and order-preserving injective maps is denoted OI. An OI-\emph{algebra} over $K$ is a (covariant) functor $\Ab$ from OI to the category of commutative, associate, unital $K$-algebras. 
\end{defn}

The category FI (resp.\ OI) is equivalent to its skeleton consisting of intervals $[n] = \{1,2,\ldots,n\}$ with $n \in \N_0$, where $[0] = \emptyset$, and injective (resp.\ order-preserving injective) maps $[m] \to [n]$. To define a functor from FI or OI to another category it is enough to define it on the corresponding skeleton. We alway use this convention. To simplify notation we often write $n$ for $[n]$ and 
$\hom_{\FI}(m,n)$  (resp. $\hom_{\OI}(m,n)$) for the set of morphisms from $[m]$ to $[n]$. Similarly, we typically use $\Ab_n$ instead of $\Ab (n)$. If $\Ab_0$ is a commutative ring then one can consider $\Ab$ as an algebra over $\Ab_0$. Thus, it is harmless to assume that $\Ab_0 = K$ as was originally required in \cite[Definition 2.4]{NR19}. All of our examples  will satisfy this condition. 

Any $\FI$-algebra $\Ab$ determines a sequence of $K$-algebras $(\Ab_n)_{n \in \N_0}$. Note that any FI-algebra is also an OI-algebra. 

The most important examples for us are the following Noetherian polynomial FI- and OI-algebras. 

\begin{defn}[{\cite[Definition 2.17]{NR19}}]
   \label{def:pol FI-algebra}
Fix any integer $c \ge 0$. 
   
(i)
Define an $\OI$-algebra $\Pb^{\OI,c}$ over $K$ by letting 
\[
\Pb^{\OI,c}_n=K\begin{bmatrix}
x_{1,1}  & \cdots & x_{1,n}\\
\vdots & \ddots & \vdots \\
x_{c,1} & \cdots & x_{c,n}
\end{bmatrix}
\]
be a polynomial ring over $K$ in $c n$ variables $x_{i, j}$ and,  for each $\epsilon\in\hom_{\OI}(m,n)$, defining a $K$-algebra  homomorphism  $\Pb^{\OI,c}(\epsilon) \colon \Pb_m\to\Pb_n$ via $x_{i,j}\mapsto x_{i,\epsilon(j)}$.

(ii) Ignoring order, we analogously define an FI-algebra $\Pb^{\FI,c}$ over $K$. 
\end{defn}

Observe that $\Pb^{\OI,c}_n=\Pb^{\FI,c}_n$ as $K$-algebras for any $n \ge 0$. 

We say that an FI-algbra $\Ab$ (resp. OI-algebra) is \emph{graded} if each $\Ab_n = \oplus_{j \in \Z} [\Ab_n]_k$ is a graded $K$-algebra and every homomorphism $\Ab (\epsilon) \colon \Ab_m \to \Ab_n$ preserves degrees.  Assigning in \Cref{def:pol FI-algebra} each variable degree one, the algebras 
$\Pb^{\FI,c}$ and $\Pb^{\OI,c}$ become graded $\FI$- and $\OI$-algebras, respectively. We always use this grading on these algebras. 

We are interested in modules over an FI-algebra. 

\begin{defn}[{\cite[Definition 3.1]{NR19}}]
  \label{def:FI-module} 
(i)
An OI-\emph{module} $\Mb$ over an OI-algebra $\Ab$ is a (covariant) functor from OI to the category of $K$-modules such that
 each $\Mb_n$ is an $\Ab_n$-module and, 
 for each $a\in\Ab_m$ and $\epsilon\in\hom_{\OI} (m,n)$, one has a commutative diagram
\begin{equation}
\label{diagram:oimodule}
\xy\xymatrixrowsep{10mm}\xymatrixcolsep{10mm}
\xymatrix {
\Mb_m\ar[d]_{a\cdot}\ar[r]^{\Mb(\epsilon)} & \Mb_n\ar[d]^{\Ab(\epsilon)(a)\cdot}\\
\Mb_m\ar[r]^{\Mb(\epsilon)} & \Mb_n
}
\endxy
\end{equation}
where the vertical maps are multiplication by the indicated elements.
We  refer to $\Mb$ also as an $\Ab$-\emph{module}. 

A morphism of OI-modules $\Mb, \Nb$ over $\Ab$ is a natural transformation $\phi \colon \Mb \to \Nb$ such that every map $\phi_n  \colon \Mb_n \to \Nb_n$ is an $\Ab_n$-module homomorphism.

(ii) Ignoring order, an FI-module over an FI-algebra and FI-module morphisms are defined analogously. 
\end{defn} 

Observe that any FI-algebra $\Ab$ is an FI-module over itself. An \emph{ideal} of $\Ab$ is a submodule of $\Ab$. 

We say that an $\Ab$-module $\Mb$ is graded if $\Ab$ is graded and each $\Mb_n$ is a graded $\Ab_n$-module. A graded morphism $\phi \colon \Mb \to \Nb$ is a map such that every $\Ab_n$-module homomorphism $\phi_n \colon \Mb_n \to \Nb_n$ preserves degrees. 

A \emph{subset} $S$ of $\Mb$, denoted $S\subseteq\Mb$, is a subset of the disjoint union $\coprod_{n\geq0}\Mb_n$. The submodule of $\Mb$ \emph{generated by a subset $S\subseteq\Mb$} is defined as the smallest submodule of $\Mb$ containing $S$ and is denoted $\la S\ra_{\Mb}$. An element of $\Mb$ is an element $f \in \Mb_n$ for some $n \ge 0$. We say that \emph{$f$ has width $n$},  denoted $w(f) = n$. We refer to $\Mb_n$ as the width $n$ component of $\Mb$. 

Any $\FI$-module determines a sequence of symmetric modules $(\Mb_n)_{n \ge 0}$. We restrict ourselves to graded $\FI$-modules over a polynomial FI-algebra $\Pb^{\FI,c}$ though analogous observations are true in greater generality. 

\begin{prop}
     \label{prop:FI gives symmetry} \mbox{ } 
     
\begin{itemize}

\item[{\rm (a)}] If $\Mb$ is a graded $\FI$-module over $\Pb^{\FI,c}$ then, for any $n \in \N_0$ and any $\sigma \in \Sym (n)$, $y \in \Mb_n$, setting 
\[
\sigma \cdot y = \Mb (\sigma) (y)
\]
defines a $K$-linear action of the symmetric group $\Sym (n)$ on $\Mb_n$ that is degree-preserving and compatible with the $\Pb^{\FI,c}_n$-module structure of $\Mb_n$. 

In particular, any $\sigma \in \Sym(n)$ acts on $\Pb^{\FI,c}_n$ as a $K$-algebra homomorphism with 
$\sigma \cdot x_{i, j} = x_{i, \sigma (j)}$. 

\item[{\rm (b)}] If $\phi \colon \Mb \to \Nb$ is a morphism of graded $\FI$-modules over $\Pb^{\FI,c}$ then, for any $n \in \N_0$, the  
$\Pb^{\FI,c}_n$-module homomorphism $\phi_n \colon \Mb_n \to \Nb_n$ is equivariant with respect to the action of $\Sym (n)$ on $\Mb_n $ and $\Nb_n$ as defined in {\rm (a)}. 

\end{itemize}

\end{prop} 

\begin{proof}

It is straightforward to check that the claims follow from the defining properties of FI-modules and their morphisms since the elements of $\Hom_{\FI} (n, n)$ are exactly the permutations on $[n]$. We leave the details to the interested reader. 
\end{proof} 

For any FI-module $\Mb$, we will always consider $\Mb_n$ as a representation of $\Sym (n)$ with the above group action.  

We use free $\OI$-modules to study analogs of free resolutions over a fixed Noetherian polynomial ring. We will consider free $\FI$-modules and discuss free resolutions of FI-modules in the next section. 

\begin{defn}[{\cite[Definition 3.17]{NR19}}]
For any OI-algebra $\Ab$ and any   integer $d \ge 0$, let $\Fb^{\OI,d}$ be the $\OI$-module over $\Ab$ defined by
\[
\Fb^{\OI,d}_n=\bigoplus_{\pi\in\hom_{\OI}(d,n)}\Ab_n e_\pi\cong (\Ab_n)^{\binom{n}{d}}   
\]
for $n \in \N_0$ and 
\[
\Fb^{\OI,d}(\epsilon) \colon \Fb^{\OI,d}_m\to\Fb^{\OI,d}_n, \;   a e_\pi\mapsto  (\Ab(\epsilon) (a)) e_{\epsilon\circ\pi},
\] 
 extended linearly, where $a \in \Ab_m$ and  $\epsilon\in\hom_{\OI}(m,n)$. 

A \emph{free} OI-module over $\Ab$ is an OI-module $\Fb$ that is  isomorphic to a direct sum 
$\bigoplus_{\lambda\in\Lambda}\Fb^{\OI,d_\lambda}$ for integers $d_\lambda\geq0$. If $|\Lambda|=n<\infty$, then $\Fb$ is said to have \emph{rank $n$}. 

To stress the dependence on $\Ab$ we also write $\Fb^{\OI,d}_{\Ab}$  for $\Fb^{\OI,d}$. 
\end{defn}

If $\Mb$ is a finitely generated OI-module over $\Pb = \Pb^{\OI,c}$ then it is Noetherian, i.e., every submodule of $\Mb$ is also finitely generated by \cite[Theorem 6.15]{NR19}. This can be used to show that there is an exact sequence of finitely generated OI-modules 
\[
\Fb^{\bullet}: \hspace*{4cm} \cdots \to \Fb^2 \to \Fb^1 \to \Fb^0 \to \Mb \to 0, 
\]
where each $\Fb^i$ is a free $\OI$-module (see \cite[Theorem 7.1]{NR19}). We refer to $\Fb^{\bullet}$ as a \emph{free $\OI$-resolution of $\Mb$} (over $\Pb$). If $\Mb$ is graded then the morphisms $\Fb^{\bullet}$ can be chosen to be graded, and the resulting resolution is called a \emph{graded free $\OI$-resolution of $\Mb$}. 

If $K$ is a field then any finitely generated $\Pb$-module $\Mb$ admits a free resolution where every free OI-module $\Fb^i$ has minimal rank compared to any other free resolution of $\Mb$. This is called a \emph{(graded) minimal free $\OI$-resolution}. It is unique up to isomorphism (see \cite[Theorem 3.10]{FN21}). 
If additionally $\Mb$ is a graded submodule of a free $\OI$-module then there is a finite algorithm that, for any integer $s \ge 0$, determines the first $s$ steps in a graded minimal free resolution of $\Mb$ by \cite[Theorem 5.4]{MN}. This algorithm relies on the theory of Gr\"obner basis for $\OI$-modules as introduced in \cite{NR19} and \cite{Nag21} (see also \cite{SS-14} for the case of OI-modules over a fixed ring). 
Our goal in the next section is to show that the first steps of an analogous free resolution of a finitely generated graded FI-module over $\Pb^{\FI,c}$ can also be determined algorithmically.


\section{Syzygies and Free Resolutions of FI-modules}
\label{section:FI-res}

The purpose of this section is to discuss how the algorithm for determining a truncated free resolution of an OI-module in \cite{MN} can be leveraged to compute   a truncated free resolution of an FI-Module.

Using the notation of \Cref{def:pol FI-algebra}, recall that $\Pb^{\OI,c}_n=\Pb^{\FI,c}_n$ as $K$-algebras for any $n \ge 0$.
 An analogous result is true in greater generality. To this end, we formalize the observation that an FI-module also defines an OI-module. Using the natural inclusion functor $I \colon \OI \to \FI$, we consider  $\OI$ as a subcategory of $\FI$. 

\begin{defn}
   \label{def:restriction to OI} 
For any FI-algebra $\Ab$ over $K$, define its \emph{restriction to \OI} as $\Ab|_\OI = \Ab\circ I$. 

For any FI-module $\Mb$ over an FI-algebra $\Ab$,  define its \emph{restriction to \OI} as $\Mb|_\OI = \Mb \circ I$. 
\end{defn}

Note that, for any $n \in \N_0$,  one has $(\Ab|_\OI)_n = \Ab_n$ as $K$-algebra and thus an equality $(\Mb|_\OI)_n = \Mb_n$ of $\Ab_n$-modules. Moreover, any morphism $\phi \colon \Mb\to\Nb$  of FI-modules, \emph{restricts} to a map $\phi|_\OI:\Mb|_\OI\to\Nb|_\OI$, which is an OI-morphism. More precisely, one has the following observation. 

\begin{prop}
   \label{prop:restriction as functor}
For any $\FI$-algebra $\Ab$, restriction gives a functor from the category of $\FI$-modules over $\Ab$ to the category of $\OI$-modules over $\Ab|_\OI$. 
\end{prop}

\begin{proof}
It is straightforward to verify the needed properties. 
\end{proof}

Recall the concept of a free FI-module.

\begin{defn}[{\cite[Definition 3.17]{NR19}}]
\label{defn:freefimodule}
Let $\Ab$ be any FI-algebra over $K$.  For any integer $d \ge 0$, define an FI-module $\Fb^{\FI,d} \colon \FI\to K\mod$ as follows: For any $n\geq0$, set
\[
\Fb^{\FI,d}_n=\bigoplus_{\pi\in\hom_\FI(d,n)}\Ab_n e_\pi\cong (\Ab_n)^{\binom{n}{d} \cdot d!}.  
\]
For each $\sigma \in\hom_\FI(m,n)$,  define $\Fb^{\FI,d}(\sigma) \colon \Fb^{\FI,d}_m\to\Fb^{\FI,d}_n$ via $ae_\pi\mapsto( \Ab(\sigma)(a))e_{\sigma\circ\pi}$. 

A \emph{free} FI-module over $\Ab$ is an FI-module $\Fb$ that is isomorphic to a direct sum $\bigoplus_{\lambda\in\Lambda}\Fb^{\FI,d_\lambda}$ for integers $d_\lambda\geq0$. 
If $|\Lambda|=r<\infty$, then $\Fb$ is said to have \emph{rank $r$}.

To stress the dependence on $\Ab$ we also write $\Fb^{\FI,d}_{\Ab}$  for $\Fb^{\FI,d}$. 
\end{defn}

Observe that free FI-modules  have the following properties. 

\begin{rem}
\label{rem:freefimodule}
Let $\Fb=\bigoplus_{i=1}^s\Fb^{\FI,d_i}$ be a free FI-module over an FI-algebra $\Ab$.
\begin{enumerate}
\item For any integer $n\geq0$,  we have
\[
\Fb_n=\bigoplus_{\substack{\pi\in\hom_\FI(d_i,n)\\1\leq i\leq s}}\Ab_ne_{\pi,i}\cong(\Ab_n)^{\sum_{i=1}^s\binom{n}{d_i}d_i!}
\]
where the second index on $e_{\pi,i}$ is for keeping track of which direct summand it lives in.
\item $\Fb$ is finitely generated as an $\Ab$-module by all $e_{\id_{[d_i]},i}$. We call these the \emph{basis elements} of $\Fb$.
\item To define an $\Ab$-linear map $\phi:\Fb\to\Nb$ where $\Nb$ is any $\Ab$-module, it suffices to specify images $\phi(e_{\id_{[d_i]},i})\in\Nb_{d_i}$ for each basis element $e_{\id_{[d_i]},i}$.
\end{enumerate}
\end{rem}

As for OI-modules, we consider subsets of and finitely generated FI-modules, and these are defined analogously. In particular, any symmetric module  generates an FI-module in the following sense. 

\begin{rem}
    \label{rem:symm mod generated FI-mod} 
 Assume $M$  is module  over a polynomial ring $\Pb^{\FI,c}_m$ with  a compatible $\Sym(m)$-action. 
    
(i) 
If $M$ is a submodule of a free $\Pb^{\FI,c}_m$-module $F$ of rank $r$ with basis $\{f_1,\ldots,f_r\}$  then denote by $\Fb$ the free 
$\Pb^{\FI,c}$-module $ \oplus_{i=1}^r \Pb^{\FI,c} e_{\id_[m], i}$. Define an $\FI$-module $\Mb = \langle M \rangle_{\Fb}$, where we consider $M$ as a subset of $\Fb_m \cong F$. Then $\Mb_m = M$ and $\Mb_n = 0$ if $n < m$. 

For example, if $c = 2$ then the ideal $I = \langle x_{1,1}^2 x_{1,2}, x_{1,1} x_{1,2}^2, x_{2,1}^4, x_{2,2}^4 \rangle \subset \Pb^{\FI,2}_2$ is a representation of $\Sym (2)$. It generates an ideal $\Ib = \langle I \rangle_{\Pb^{\FI,2}}$ of $\Pb^{\FI,2}$ with 
\[
\Ib_n = \langle x_{1,i}^2 x_{1,j}, x_{2, k}^4 \; \mid \; k \in [n], \ i \neq j \in [n]\}. 
\]

(ii) If $M$ is finitely generated then pass to a syzygy module $Q$ given by an exact sequence of finitely generated $\Pb^{\FI,c}_m$-modules
\[
0 \to Q \stackrel{\phi}{\longrightarrow} F \to M \to 0,  
\]
where $F$ is a free module. Let $\Nb$  be the FI-module generated by $\phi (Q)$ as a submodule of the free FI-module $\Fb$ with $\Fb_m \cong F$ as in (a). Then $\Fb/\Nb$ is an FI-Module with $(\Fb/\Nb)_m \cong M$. 
\end{rem}

The following result relates finitely generated free FI-modules with finitely generated free OI-modules.

\begin{prop}
\label{prop:freefibecomesfreeoi}
Fix integers $d_1,\ldots,d_s\geq0$ and let $\Fb=\bigoplus_{i=1}^s\Fb^{\FI,d_i}$ be a free \FI-module over an \FI-algebra $\Ab$. 
There is an isomorphism of $\OI$-modules
\[
\Fb|_\OI \cong  \bigoplus_{i=1}^s(\Fb^{\OI,d_i})^{\oplus d_i!} 
\]
over $\Ab|_\OI$, and thus $\Fb|_\OI$ is a free \OI-module of rank $\sum_{i=1}^s d_i!$.
\end{prop}

\begin{proof} Recall that each module $\Fb^{\OI,d_i}$ has rank one. Hence, 
the free OI-module $\bigoplus_{i=1}^s(\Fb^{\OI,d_i})^{\oplus d_i!}$ has rank $\sum_{i=1}^s d_i!$. Fix a basis 
\[
\{e_{\id_{[d_i]},i,j}\;:\; i\in[s], j\in[d_i!]\}.
\]
For each $i\in[s]$, the set $\Sym(d_i)=\hom_{\FI}(d_i,d_i)$ has $d_i!$ elements, which we denote by $\{\pi_{i,1},\ldots,\pi_{i,d_i!}\}$. Fix a basis $\{e_{\id_{[d_i]},i}\;:\;i\in[s]\}$ of $\Fb$ and define an OI-morphism \[\phi:\bigoplus_{i=1}^s(\Fb^{\OI,d_i})^{\oplus d_i!}\to\Fb|_\OI\quad\text{via}\quad
e_{\id_{[d_i]},i,j}\mapsto e_{\pi_{i,j},i}.
\]
Consider any integer $n\geq0$ and note that $(\bigoplus_{i=1}^s(\Fb^{\OI,d_i})^{\oplus d_i!})_n$ has an $(\Ab|_\OI)_n$-basis given by
\[
A=\{e_{\sigma,i,j}\;:\; i\in[s], j\in[d_i!], \sigma\in\hom_{\OI}(d_i,n)\}.
\]
Furthermore, $(\Fb|_\OI)_n$ has an $(\Ab|_\OI)_n$-basis given by
\[
B=\{e_{\tau,i}\;:\;i\in[s],\tau\in\hom_{\FI}(d_i,n)\}.
\]
It suffices to show that $\phi_n$ induces a bijection from $A$ to $B$. Note that $\phi_n(e_{\sigma,i,j})=e_{\sigma\circ\pi_{i,j},i}$ and $\sigma\circ\pi_{i,j}$ is an FI-morphism from $[d_i]$ to $[n]$, and so $\phi_n(A)\subseteq B$. Now pick any $e_{\tau,i}\in B$. There is a unique $\sigma\in\hom_{\OI}(d_i,n)$ with the same image as $\tau$. Thus $\tau=\sigma\circ\pi_{i,j}$ for some $\pi_{i,j}\in \Sym(d_i)$,  
which implies $\phi_n(e_{\sigma,i,j})=e_{\sigma\circ\pi_{i,j},i}=e_{\tau,i}$ and hence shows that $\phi_n$ surjects $A$ onto $B$. 
Since $A$ and $B$ both share the same finite cardinality of $\sum_{i=1}^s\binom{n}{d_i}d_i!$, it follows that $\phi_n$ bijects $A$ onto $B$, and we are done.
\end{proof}

\begin{rem}
The isomorphism in \Cref{prop:freefibecomesfreeoi} is not canonical. In particular, it depends on how one chooses to order the sets $\Sym(d_i)$ for $i\in[s]$.
\end{rem}

Given a submodule $\Mb$ of a finitely generated free FI-module $\Fb$ over $\Pb^{\FI,c}$, there is an induced submodule $\Mb|_\OI$ of the finitely generated free OI-module $\Fb|_\OI$ over $\Pb^{\OI,c}$. Thus, one could immediately apply the results of the previous chapter to compute a free OI-resolution of $\Mb|_\OI$. However, such a resolution will only be compatible with the OI-symmetries of $\Mb$. In order to resolve $\Mb$ in a way compatible with all its FI-symmetries, one desires a free FI-resolution
\begin{equation*}
   \label{eq:freefires}
\Fb^\bullet: \quad \cdots\to\Fb^2\to\Fb^1\to\Fb^0\to\Mb\to0, 
\end{equation*}
where each $\Fb^i$ is a finitely generated free FI-module over $\Pb^{\FI,c}$ and the differentials are FI-morphisms. By \cite[Theorem 7.1]{NR19}, such a free FI-resolution of $\Mb$ always exists. If $\Mb$ is graded and $K$ is a field, the arguments in \cite{FN21} can be used to show that $\Mb$ admits a graded such resolution $\Fb^\bullet$ that is \emph{minimal}, that is, each free module $\Fb^i$ in $\Fb^\bullet$ has minimal rank compared to the free module in homological degree $i$ of any other graded free $\FI$-resolution of $\Mb$. Such a graded free FI-resolution is unique up to isomorphism. 

The existence argument for  $\Fb^\bullet$ given in \cite{NR19} is not constructive. In what follows, we present the first algorithmic method for computing a free FI-resolution of $\Mb$ in which each free $\Pb^{\FI,c}$-module is finitely generated, provided $K$ is a field. 

We need some more preparatory observations.

\begin{lem}
\label{lem:fioiker}
Let $\phi:\Mb\to\Nb$ be a map of \FI-modules over an \FI-algebra $\Ab$. Suppose $\ker(\phi|_\OI)$ is generated as an $\Ab|_\OI$-module by a subset $B\subseteq\Mb|_\OI$. Then $\ker(\phi)$ is generated as an $\Ab$-module by $B\subseteq\Mb$.
\end{lem}

\begin{proof}
First note that $B\subseteq\ker(\phi)$ since for all $b\in B$, we have $\phi(b)=\phi|_\OI(b)=0$. Now pick any $f\in\ker(\phi)$. Then $\phi|_\OI(f)=\phi(f)=0$, and so $f\in\ker(\phi|_\OI)$. Thus,  there is an expression $f=\sum a_i \tau_i (q_i)$ with elements 
$a_i\in(\Ab|_\OI)_{w(f)}$,  
$q_i \in B$ and $\tau_i \in \Hom_{\OI} (w(q_i), w(f))$. 
The result follows by noticing that $(\Ab|_\OI)_{w(f)}=\Ab_{w(f)}$ and that $\tau_i$ is also an element of $\Hom_{\FI} (w(q_i), w(f))$.
\end{proof}

For the remainder of this section, we assume that $K$ is a field. This allows us  to apply results in \cite{MN}, where a special case of the following result has been established.  

\begin{prop}
\label{lem:oikergen}
Let $\phi \colon \Gb\to\Fb$ be a map of finitely generated free \OI-modules over $\Pb=\Pb^{\OI,c}$. Then a finite generating set of 
$\ker(\phi)$ can be computed in finite time.
\end{prop}

\begin{proof}
Fix a basis $\{\varepsilon_{\id_{[d_i]},i}\;:\;i\in[s]\}$  of $\Gb = \bigoplus_{i=1}^s F^{\OI, d_i}$ with integers $d_1,\ldots,d_s\geq0$. For each $i\in[s]$, set $b_i=\phi(\varepsilon_{\id_{[d_i]},i})\in\Fb_{d_i}$. Choose a monomial order $<$ on $\Fb$ and let $B=\{b_1,\ldots,b_s\} \subset\Fb$. One checks in finite time whether $B$ is a Gröbner basis of $\la B\ra_{\Fb}$ with respect to $<$  by using the OI-Buchberger's criterion \cite[Theorem 3.16]{MN}. If $B$ is a Gröbner basis then the result follows immediately by applying the OI-Schreyer's theorem \cite[Theorem 4.6]{MN} to $\phi$.

Suppose now that $B$ is not a Gröbner basis. We will run the OI-Buchberger's algorithm on $B$, see \cite[Algorithm 3.17]{MN}.  Set $B_0=B$ and let $\Fb(\sigma_1)(b_{i_1}),\Fb(\tau_1)(b_{j_1})$ be elements of $\langle B_0 \rangle_{\Fb}$ in some width $m_1$ whose S-polynomial has a nonzero remainder $b_{s+1} \in\Fb_{m_1}$ modulo $B_0$. This gives an expression
\begin{equation}
\label{eq:firstdivision}
S(\Fb(\sigma_1)(b_{i_1}),\Fb(\tau_1)(b_{j_1}))=\sum_\l a_{1,\l} q_{1,\l} +b_{s+1}
\end{equation}
with $a_{1,\l}\in\Pb_{m_1}$ and $q_{1,\l}\in (\langle B_0 \rangle_{\Fb})_{m_1}$. 
Now set $B_1=B_0\cup\{b_{s+1}\}$ and repeat this process until the algorithm terminates to obtain sets $B_0\subset\cdots\subset B_t$, 
where $B_t = \{b_1,\ldots,b_s,b_{s+1},\ldots,b_{s+t}\}$ is a Gröbner basis of $\langle B \rangle_{\Fb}$ with respect to $<$ and we get expressions
\begin{equation}
\label{eq:expressions}
S(\Fb(\sigma_v)(b_{i_v}),\Fb(\tau_v)(b_{j_v}))=\sum_\l a_{v,\l} q_{v,\l} +b_{s+v}
\end{equation}
with $a_{v,\l}\in\Pb_{m_v}$ and $q_{v,\l}\in (\langle B_{v-1} \rangle_{\Fb})_{m_v}$ for each $v\in[t]$.

Since $B_t$ is a Gr\"obner bases, \cite[Theorem 3.16]{MN} yields that, for any $b_i, b_j \in B_t$ and any $\sigma \in \Hom_{\OI} (w(b_i), m)$ and $\tau \in \Hom_{\OI} (w(b_j), m)$ with $m = | \im \sigma \cup \im \tau |$, any nonzero S-polynomial of $\Fb(\sigma)(b_i)$ and $\Fb(\tau)(b_j)$ has a remainder of zero modulo $B_t$, that is, there is a relation of the form 
\[
S(\Fb(\sigma)(b_i),\Fb(\tau)(b_j))=\sum_\l a_{i,j,\l}^{\sigma,\tau}\Fb(\pi_{i,j,\l}^{\sigma,\tau})(b_{k_{i,j,\l}^{\sigma,\tau}})
\]
with $a_{i,j,\l}^{\sigma,\tau}\in\Pb_m$, $b_{k_{i,j,\l}^{\sigma,\tau}} \in G_t$  and $\Fb(\pi_{i,j,\l}^{\sigma,\tau})(b_{k_{i,j,\l}^{\sigma,\tau}})\in \Fb_m$. Each such relation defines an element
\begin{equation}
\label{eq:syz}
s_{i,j}^{\sigma,\tau}=m_{i,j}^{\sigma,\tau}\varepsilon_{\sigma,i}-m_{j,i}^{\tau,\sigma}\varepsilon_{\tau,j}-\sum_\l a_{i,j,\l}^{\sigma,\tau}\varepsilon_{\pi_{i,j,\l}^{\sigma,\tau},k_{i,j,\l}^{\sigma,\tau}}\in\Gb_m
\end{equation}
with $m_{i,j}^{\sigma,\tau}, m_{j,i}^{\tau,\sigma} \in \Pb_m$. 

Now set $d_{s+v}=m_v=w(b_{s+v})$ for $v\in[t]$ and let $\widehat{\Gb}$ be the free $\Pb$-module with basis $\{\varepsilon_{\id_{[d_i]},i}\;:\;i\in[s+t]\}$. Observe that this basis contains the above basis of $\Gb$.  Define a $\Pb$-linear map
\[
\widehat{\phi} \colon \widehat{\Gb}\to\Fb\quad\text{via}\quad \varepsilon_{\id_{[d_i]},i}\mapsto b_i
\]
and equip $\widehat{\Gb}$ with the Schreyer order $<_{B_t}$ induced by $<$ and $B_t$. By construction, each of the above elements $s_{i,j}^{\sigma,\tau}$ is in the kernel of $\widehat{\phi}$, that is, $s_{i,j}^{\sigma,\tau}$ is a syzygy of $\langle B \rangle_{\Fb}$. 
 In fact, by the OI-Schreyer's theorem \cite[Theorem 4.6]{MN}, $\ker(\widehat{\phi})$ is generated by these syzygies $s_{i,j}^{\sigma,\tau}$.  
Denote by  $\mcal{S}$ this finite generating set of $\ker(\widehat{\phi})$. 

Notice that the expressions from \eqref{eq:expressions} give rise to syzygies $s_{i_v,j_v}^{\sigma_v,\tau_v}\in\ker(\widehat{\phi})$ with the following special form:
\begin{equation}
\label{eq:specialform}
s_{i_v,j_v}^{\sigma_v,\tau_v}=\varepsilon_{\id_{[d_{s+v}]},s+v}+\sum_{\ell} b_{v, \ell} \varepsilon_{\rho_{v, \ell},  u_{v, \ell}}
\end{equation}
with $b_{v, \ell} \in \Pb_{d_{s+v}}$ and $u_{v, \ell} <s+v$  for each $\ell$.

Denote by $\Mb$ be the submodule of $\ker(\widehat{\phi})$ generated by the syzygies 
$s_{i_v,j_v}^{\sigma_v,\tau_v}$ with $v\in[t]$. Passing to the quotient $\widehat{\Gb}/\Mb$, there is an induced map 
$\psi \colon \widehat{\Gb}/\Mb\to\Fb$ given by $g+\Mb \mapsto \widehat{\phi}(g)$. If $\alpha \colon \widehat{\Gb}\to\widehat{\Gb}/\Mb$ denotes the canonical projection map, then $\alpha(\mcal{S})$ is a generating set of $\ker(\psi)$. 
Using \eqref{eq:specialform}, we see that there is an isomorphism $\Phi \colon \Gb\to\widehat{\Gb}/\Mb$ given by $\varepsilon_{\id_{[d_i]},i}\mapsto\alpha (\varepsilon_{\id_{[d_i]},i})$ for each $i\in[s]$. Since $\phi=\psi\circ\Phi$, a finite generating set for $\ker(\phi)$ is given by $\Phi^{-1}(\alpha (\mcal{S}))$, as desired.
\end{proof}

Putting everything together, we obtain an algorithm for computing truncated free FI-resolutions over a Noetherian polynomial FI-algebra.

\begin{proc}[Method for Computing a Free FI-Resolution]
\label{proc:freefires} \mbox{ }

Input: A finite subset $B$ of a finitely generated free FI-module $\Fb = \bigoplus_{i = 1}^r \Fb^{\FI, d_i}_{\Pb} $ over $\Pb=\Pb^{\FI,c}$, $p \ge 1$ an integer.

Output: An exact sequence of finitely generated OI-modules over $\Pb$, 
\[
 \Fb^p \to \Fb^{p-1} \to \cdots \to  \Fb^1 \to \Mb \to 0
\]
with $\Mb = \langle B \rangle_{\Fb}$ and free OI-modules $\Fb^i$.

\begin{enumerate}
\item Define a morphism of free $\Pb$-modules  
\[
\phi_1 \colon \Fb^1 = \bigoplus_{i = 1}^s \Fb^{\FI, w (b_i)}_{\Pb} \to \Fb^0 = \Fb \;  \text{ by }  \phi_1 (\eps_{ \id_{w (b_i)}, i}) \mapsto b_i, 
\]
where $B_0 = B = \{b_1,\ldots,b_s\}$ and $\{\eps_{ \id_{w (b_i)}, i} \; \mid \; i \in [s]\}$ is a basis of $\Fb^1$. Notice that the image of $\phi_1$ is $\Mb$. 

\item Apply \Cref{prop:freefibecomesfreeoi} to obtain a map $\phi_1|_\OI:\Fb^1|_\OI\to\Fb^{0}|_\OI$ of finitely generated free OI-modules over $\Pb^{\OI,c}$.
\item Use \Cref{lem:oikergen} to compute a finite generating set $B_1$ of $\ker(\phi_1|_\OI)$. By \Cref{lem:fioiker}, $B_1$ generates $\ker(\phi_1)$ as an FI-module over $\Pb^{\FI,c}$.

\item Repeat steps (i)-(iii) $p-1$ times to obtain the beginning of a free FI-resolution of $\Mb$, 
\[
 \Fb^p \to \Fb^{p-1} \to \cdots \to  \Fb^1 \to \Mb \to 0. 
\]

\end{enumerate}
\end{proc}

\begin{rem}
   \label{rem:graded FI-res}
If $\Fb$ is a graded $\FI$-module and the elements of $B$ are homogeneous then $\Mb$ is a graded $\Pb$-module. Choosing the gradings of the free FI-modules $\Fb^i$ suitably, \Cref{proc:freefires} gives the beginning of a graded free FI-resolution. For example, one sets $\Fb^1 = \bigoplus_{i = 1}^s \Fb^{\FI, w (b_i)}_{\Pb} (-t_i)$ if $t_i$ is the degree of $b_i$, where,  for any graded module $N$ over a graded ring $A$, the degree shifted module $N (t)$ is the $A$-module with the same $A$-module structure as $N$ but grading given by $[N(t)]_j = [N]_{t +j}$ for $j \in \Z$. With these adjustments, the morphisms $\phi \colon \Fb^i \to \Fb^{i-1}$ are then graded morphisms. 

Furthermore, as mentioned above,  $\Mb$ admits a graded minimal free $\FI$-resolution. If desired such a minimal resolution can be obtained by adding to \Cref{proc:freefires}  additional steps that modify, beginning with $\phi_1$, the obtained morphisms $\phi_i$ to minimal morphisms.  The pruning is analogous to the case of $\OI$-morphisms, which is described in detail in the proof of  \cite[Lemma 3.5]{FN21}. 
\end{rem}

In the graded case, any graded free FI-resolution of $\Mb$ over $\Pb$ gives simultaneously, for any $n \ge 0$, an  equivariant graded free resolution of $\Mb_n$ over $\Pb_n$ by \Cref{prop:FI gives symmetry}.  Such a resolution of $\Mb_n$ can be modified to obtain an equivariant graded minimal free resolution by \Cref{prop:prune res} (in the non-modular case).


\section{Colimits of OI-modules}
\label{section:Inc-resolutions}

Denote the monoid of increasing functions  $\N \to \N$ by $\Incn$ or simply $\Inc$. 
We introduce Inc-modules over Inc-algebras and show how important examples arise as colimits of OI-modules. 
Some steps in this direction were taken in \cite{HS, NR19}. However, in these sources the framework was neither fully developed nor utilized for investigating modules other than ideals. We use our set-up to establish algorithms for  computing $\Inc$-equivariant  free resolutions of finitely generated Inc-modules over a Noetherian polynomial Inc-algebra. 

Denote by $K$ any commutative ring.

\begin{defn}
   \label{def:Inc-algebra}
(i)
An \emph{$\Inc$-algebra} over $K$ is a $K$-algebra $\Ac$ that admits an action $\Inc \times \Ac \to \Ac, \tau \times a \mapsto \tau \cdot a$,  such that $\tau(-) \colon \Ac \to \Ac$,  $\tau (a) = \tau \cdot a$, is a $K$-algebra homomorphism for any $\tau \in \Inc$.  

(ii) An \emph{$\Inc$-algebra homomorphism} is an Inc-equivariant $K$-algebra homomorphism $\phi \colon \Ac \to \Bc$ of Inc-algebras, that is, it satisfies $\phi (\tau \cdot a) = \tau \cdot \phi (a)$ for any $\tau \in \Inc$ and $a \in \Ac$. 
\end{defn}

If $\Ac$ is a monoid ring over $K$  one may work with the skew-monoid ring $\Ac * \Inc$ as introduced in \cite{HS} instead of considering $\Ac$ an an $\Inc$-algebra. 

\begin{ex}
    \label{exa:pol as Inc-algebr}
(i)
The polynomial ring $K[X_{c \times \N}]$ is an Inc-algebra with an action \\
$\Inc \times K[X_{c \times \N}] \to K[X_{c \times \N}] $ that is induced by 
\[
\tau \cdot x_{i, j} := x_{i, \tau (j)}. 
\]

(ii)
The projection $K[X_{c \times \N}] \to K[X_{(c-1) \times \N}]$ induced by 
\[
x_{i, j} \mapsto \begin{cases}
x_{i, j} & \text{ if } i \in [c-1]; \\
0 & \text{ if } i = c
\end{cases}
\]
is an Inc-algebra homomorphism. 
\end{ex} 

Our first goal is to show how Inc-algebras arise as colimits of OI-algebras. 
Recall that we often  write  $\Hom_{\OI} (m, n)$  for  $\hom_{\OI}([m],[n])$ and  $\id_n$ for the identity map $\id_{[n]}$ on $[n]$. 

For any integer $d \ge 0$,  we consider the set of order-preserving maps from $[d]$ to $\N$,  
\[
\Hom_{\OI} (d, \N) = \{ \pi \colon d \to \N \; \mid \; \pi(1) < \pi (2) < \cdots < \pi (d) \}.
\]  
Composition provides an action of $\Inc$ on this set. 

\begin{prop}

For any integer $d \ge 0$, the map
\[
\Inc \times \Hom_{\OI} (d, \N) \to \Hom_{\OI} (d, \N), \; \tau \times \pi \mapsto \tau \circ \pi, 
\]
gives an action of $\Inc$ on $\Hom_{\OI} (d, \N)$. 
\end{prop}

\begin{proof}
Since $\tau$ and  $\pi$ preserve the order on $\N$ so does $\tau \circ \pi$. Thus, the given map is well-defined. It satisfies the properties of an action. 
\end{proof} 

\begin{rem}
    \label{rem:hom set as colimit}
Consider the direct system $( \Hom_{\OI} (d, n), \iota_{m, n})$ with maps  
\begin{align}
\label{eq:iota}
\iota_{m, n}\colon [m] \to [n], \ j \mapsto j. 
\end{align} 
Its colimit is the set $\Hom_{\OI} (d, \N)$, where the class of $\pi \in \Hom_{\OI} (d, n)$ is the map $\tilde{\pi} \colon d \to \N$ with 
$\tilde{\pi}(j) = \pi (j)$ for $j \in [d]$. 
\end{rem}

In view of the above remark, the following definition is rather natural. 

\begin{defn}
   \label{def:action of Inc on OI-maps}
For any $\tau \in \Inc$ and any $\pi \in \hom_{\OI} (d, n)$, define a map
\[
\tau \cdot \pi \colon [d] \to [\tau (\pi (n))], \; j \mapsto \tau (\pi (j)). 
\] 
It is in $\Hom_\OI (d, \tau (\pi (n)))$. 
\end{defn}

Using also the maps $\iota_{m, n}$ as defined in \Cref{eq:iota}, we get the following result. For ease of notation, given an OI-algebra $\Ab$ and any $\sigma \in \hom_{\OI} (m, n)$, we use $\sigma_*$ to denote the induced map $\Ab(\sigma) \colon \Ab_m \to \Ab_n$. 

\begin{prop}
    \label{prop:colomit of OI-alg}
Given any $\OI$-algebra $\Ab$,
using the direct system $(\Ab_n, \iota^*_{m, n})$,  its colimit $\Ac = {\displaystyle \lim_{\longrightarrow}} \, \Ab_n$ 
is an $\Inc$-algebra whose $\Inc$-action is defined by
\[
\tau \cdot [a] = [(\tau \cdot \id_n)_* (a)], 
\]
where $[a]$ denotes the class of $a \in \Ab_n$ and $\tau \in \Inc$. 
\end{prop} 

\begin{proof}
This follows from the arguments for \cite[Proposition 2.10]{NR19}. 
\end{proof}

Throughout this paper, we always use the above direct system and Inc-action for colimits of an OI-algebra. 

\begin{ex}
  \label{exa:poly Inc-algebras}
(i)
Given any polynomial OI-algebra $\Ab = X^{\OI, d}$ with $d \ge 0$ (see \cite[Definition 2.17]{NR19}), its colimit ${\displaystyle \lim_{\longrightarrow}} \, \Ab = {\displaystyle \lim_{\longrightarrow}}  \, \Ab_n$ is the  polynomial $K$-algebra $\Ac = K[x_{\pi} \; \mid \; \pi \in \Hom_{\OI} (d, \N)]$ whose $\Inc$-action is induced by 
\[
\tau \cdot x_{\pi} = x_{\tau \circ \pi}. 
\]

(ii)
For any integer $c \ge 1$, the colimit of $\Pb = \Pb^{\OI, c}$ is isomorphic to the Inc-algebra defined in \Cref{exa:pol as Inc-algebr}(i), where the class $[f]$ of any $f \in \Pb_n$ is mapped onto $f$, considered as a polynomial in $K[X_{c \times \N}]$. 
Indeed, this follows by identifying any variable $x_{i, \pi }$ with $x_{i, \pi  (1)}$ for any  $\pi  \in \Hom_{\OI} (1, \N)$ and $i \in [c]$. 
\end{ex}

Taking the Inc-action into account, there is a natural concept of finite generation. 

\begin{defn}
   \label{def:fg Inc-algebra} 
An Inc-algebra $\Ac$ is said to be \emph{finitely generated} 
if there exists
a finite subset $G \subset  \Ac$
which is not contained in any proper Inc-subalgebra of $\Ac$.
\end{defn}

\begin{ex}
  \label{exa:fg Inc algebra}
The polynomial ring $K[x_{\pi } \; \mid \; \pi  \in \Hom_{\OI} (d, \N)]$ considered in \Cref{exa:poly Inc-algebras} is not a finitely generated $K$-algebra. However, as an Inc-algebra it is generated by one variable, namely $x_{\widetilde{\id_d}}$. See \Cref{rem:hom set as colimit} for notation. 
\end{ex}

We are ready to introduce our main objects of interest in this section. 

\begin{defn}
   \label{def:Inc-module} 
Let $\Ac$ be an Inc-algebra over $K$. The category of \emph{Inc-modules over} $\Ac$ is denoted by $\Inc$-$\Mod (\Ac)$. An object of  $\Inc$-$\Mod (\Ac)$ is any $\Ac$-module $\Mcc$ that admits an action $\Inc \times \Mcc \to \Mcc, \ \tau \times q \mapsto \tau \cdot q$, such that 
\begin{itemize}

\item[(i)] $\tau \cdot (q_1 + q_2) = \tau \cdot q_1 + \tau \cdot q_2$ \;  if $\tau \in \Inc$ and $q_1, q_2 \in \Mcc$; and

\item[(ii)] $\tau \cdot (a q) = (\tau \cdot a) (\tau \cdot q)$ \; if $\tau \in \Inc$, $a \in \Ac$ and $q \in \Mcc$;
\end{itemize}

A morphism of  $\Inc$-$\Mod (\Ac)$ is any Inc-equivariant homomorphism $\phi \colon \Mcc \to \Nc$ of $\Ac$-modules, that is, it satisfies $\phi (\tau \cdot q) = \tau \cdot \phi (q)$ for any $\tau \in \Inc$ and $q \in \Mcc$. 
\end{defn}

Notice that any Inc-algebra $\Ac$ is an Inc-module over itself. 

\begin{rem}
One checks for any $\tau \in \Inc$, $a \in \Ac$ and $q \in \Mcc$ that 
\[
\tau \cdot (a \phi (q)) = (\tau \cdot a) (\tau \cdot \phi (q)). 
\]
It follows that the image of $\phi$ is an Inc-submodule of $\Nc$, and so $\Inc$-$\Mod (\Ac)$ is an abelian category. 
\end{rem} 

We now introduce a family of Inc-modules that will be important for us. 

\begin{defn} 
   \label{def:basic free Inc-module}
For any   integer $d\geq0$, define an Inc-module $\Fc^{\Inc,d} = \Fc^{\Inc,d}_{\Ac}$ over an Inc-algebra $\Ac$ as the free $\Ac$-module 
\[
\Fc^{\Inc,d} =\bigoplus_{\pi\in\Hom_{\OI} (d,\N)}\Ac \eps_\pi   
\]
with basis $\{ \eps_{\pi} \; \mid \; \pi \in \Hom_{\OI} (d,\N) \}$ 
together with the action $\Inc \times \Fc^{\Inc,d} \to \Fc^{\Inc,d}$ induced by
\[
\tau \cdot (a e_{\pi}) = (\tau \cdot a) e_{\tau \circ \pi}. 
\]
One checks that this defines indeed an Inc-module over $\Ac$. 
\end{defn} 

Note that $\Fc^{\Inc,0}_{\Ac} \cong \Ac$ because $[0]$ is the empty set. 

In order to study properties of the above and other Inc-modules, the following concepts will be useful. 

\begin{defn}
   \label{def:fg Inc-module} 

(i) For any subset $G$ of an Inc-module $\Mcc$, the Inc-submodule of $\Mcc$ \emph{generated by $G$} is the smallest Inc-submodule of $\Mcc$ that contains $G$. It is denoted by $\langle G \rangle$ or, more precisely, by $\langle G \rangle_{\Mcc}$. 

(ii) 
An Inc-module $\Mcc$ is said to be \emph{finitely generated} if there is a finite subset $G \subset \Mcc$ such that $\langle G \rangle_{\Mcc} = \Mcc$. 

(iii) An Inc-module is \emph{Noetherian} if every Inc-submodule is finitely generarated. 
\end{defn}

Explicitly, the elements of  $\langle G \rangle_{\Mcc} \subseteq \Mcc$ are finite sums 
\[
\sum_i a_i  (\tau_i \cdot g_i) 
\]
with suitable $g_i \in G$, $\tau_i \in \Inc$ and $a_i \in \Ac$. 

Observe that each module $\Fc^{\Inc,d}$ is generated by one element as an Inc-module, namely by  $\eps_{\widetilde{\id_d}}$, where $\widetilde{\id_d}\colon [d] \to \N, j \mapsto j$, see \Cref{rem:hom set as colimit}. Moreover, any morphism $\Fc^{\Inc,d}  \to \Mcc$ of Inc-modules is determined by the image of $\eps_{\widetilde{\id_d}}$. This justifies the following terminology. 

\begin{defn} 
    \label{def:free Inc-module}
A \emph{free} Inc-module over $\Ac$ is any Inc-module $\Fc$ that is isomorphic to a direct sum 
$\bigoplus_{\lambda\in\Lambda}\Fc^{\Inc,d_\lambda}_{\Ac} $ for integers $d_\lambda\geq0$. If $|\Lambda|=r<\infty$, then $\Fc$ is said to have \emph{rank $r$}. 
\end{defn} 

As in the classical case, finite generation is a consequence of the existence of a suitable surjection. 

\begin{prop}
    \label{prop:fg Inc-module}
If  $\Mcc$ is an Inc-module that admits a surjection
\[
\bigoplus_{i = 1}^k \Fc^{\Inc,d_i}  \to \Mcc
\]
with suitable integers $d_i \ge 0$ then $\Mcc$ is a finitely generated Inc-module.
\end{prop}

\begin{proof}
The images of the elements $e_{\widetilde{\id_{d_i}}}$ generate $\Mcc$. 
\end{proof}

Below, we will consider a case, where the converse is true as well. 

\begin{prop}
    \label{prop:colomit of OI-mod}
For any $\OI$-algebra $\Ab$ and any $\OI$-module $\Mb$ over $\Ab$,  taking the colimit $\Mcc = {\displaystyle \lim_{\longrightarrow}} \, \Mb$ 
using the direct system $(\Mb_n, \Mb(\iota_{m, n}))$ with maps  
\begin{align*}
\label{eq:iota}
\iota_{m, n}\colon [m] \to [n], \ j \mapsto j,
\end{align*}
is an exact functor from the category of $\OI$-modules over $\Ab$ into the category of Inc-modules over $\Ac = {\displaystyle \lim_{\longrightarrow}} \, \Ab$. The scalar multiplication on $\Mcc$ is defined by 
\[
[a] \cdot [q] = [ \Ab (\iota_{m,k})(a) \Mb(\iota_{n, k}) (q) ], 
\]
where $a \in \Ab_m$, $q \in \Mb_n$ and $k = \max \{m, n\}$, and the Inc-action is defined by 
\[
\tau \cdot [q] = [ \Mb(\tau \cdot \id_n) (q)], 
\]
where $\tau \in \Inc$, $q \in \Mb_n$ and $\tau \cdot \id_n \colon [m] \to [\tau (n)], j \mapsto \tau (j)$, see \Cref{def:action of Inc on OI-maps}. 

Any $\OI$-morphism $\phi \colon \Mb \to \Nb$ induces an Inc-module morphism $\overline{\phi} \colon \Mcc = {\displaystyle \lim_{\longrightarrow}} \, \Mb \to \Nc = {\displaystyle \lim_{\longrightarrow}} \, \Nb, [q] \mapsto [\phi (q)]$. 
\end{prop} 

\begin{proof} The arguments for \cite[Lemma 3.5 and Lemma 3.6]{NR19} show that $\Mcc$ is indeed an Inc-module over $\Ac$ with the stated scalar multiplication and Inc-action. 
The induced map $\overline{\phi} \colon \Mcc \to \Nc$ is Inc-equivariant because 
\[
[\phi (\Mb ( \tau \cdot \id_n)(q)] = [\Nb (\tau \cdot \id_n)(\phi (q)]. 
\] 
Exactness of the colimit in our situation is a standard fact (see, e.g., \cite[Theorem 2.6.15]{W}). 
\end{proof}

The next result identifies colimits of free OI-modules and implies that these are free Inc-modules. Recall that any map $\pi \in \hom_{\OI} (m, n)$ defines a map $\tilde{\pi} \in \hom_{\OI} (m, \N)$, as introduced in \Cref{rem:hom set as colimit}. 

\begin{prop} 
    \label{prop:colimit of free OI-module} 
For any $\OI$-algebra $\Ab$, one has 
\[
\lim_{\longrightarrow} \Fb^{\OI,d}_{\Ab} = \Fc^{\Inc,d}_{\Ac}
\]
with $\Ac \cong {\displaystyle \lim_{\longrightarrow}} \, \Ab$. 
\end{prop}

\begin{proof}
It follows by a  routine argument that the map 
\[
\lim_{\longrightarrow} \Fb^{\OI,d}_{\Ab}  \to \Fc^{\Inc,d}_{\Ac} \quad \text{ defined by }   [a e_{\pi}]  \mapsto [a] \eps_{\tilde{\pi}} 
\]
is an isomorphism. We leave the details to the interested reader. 
\end{proof}

Since any finitely generated OI-module $\Mb$ admits a surjection $\Fb \to \Mb$ with a finitely generated free OI-module $\Fb$ by  \cite[Proposition 3.18]{NR19}, we immediately obtain a conceptual proof for the following result. 

\begin{cor}[{\cite[Proposition 3.14]{NR19} }]
If $\Mb$ is a finitely generated $\OI$-module over $\Ab$ then its colimit $\Mcc = {\displaystyle \lim_{\longrightarrow}} \, \Mb$ is a finitely generated $\Inc$-module over $\Ac =  {\displaystyle \lim_{\longrightarrow}} \, \Ab$.
\end{cor} 

Noetherianity transfers as well.  

\begin{cor}
   \label{cor:limit is noeth}
If $\Mb$ is a Noetherian $\OI$-module over $\Ab$ then its colimit $\Mcc =  {\displaystyle \lim_{\longrightarrow}} \, \Mb$ is a Noetherian $\Inc$-module over $\Ac = {\displaystyle \lim_{\longrightarrow}} \, \Ab$.
\end{cor} 

\begin{proof}
Combine \Cref{prop:colomit of OI-mod}
and \cite[Theorem 4.6]{NR19}. 
\end{proof}

For the remainder of this section, we fix an integer $c \ge 1$ and consider Inc-modules over the polynomial Inc-algebra $\Pc = K[X_{c \times \N}] \cong {\displaystyle \lim_{\longrightarrow}} \, \Pb$ with  $\Pb =  (\Pb^{\OI, c})$, as introduced in \Cref{exa:pol as Inc-algebr}.  
Main results in \cite{C} or \cite{HS} show that any Inc-ideal of $\Pc$ is finitely generated, that is, $\Pc$ is Noetherian. This is true more generally for submodules of finitely generated free Inc-modules.

\begin{cor}
     \label{cor:fg free Inc-mod is noeth}
If $K$ is a Noetherian ring then any finitely generated free Inc-module over $\Pc = K[X_{c \times \N}]$ is Noetherian.
\end{cor}

\begin{proof} 
Consider any finitely generated free Inc-module $\Fc$ over $\Pc$.  Thus, $\Fc$ is isomorphic to $\bigoplus_{i = 1}^r \Fc^{\Inc, d_i}_{\Pc} $ with integers $d_i \ge 0$. \Cref{prop:colimit of free OI-module} gives that $\Fc =  {\displaystyle \lim_{\longrightarrow}} \, \Fb$ with $\Fb = \bigoplus_{i = 1}^r \Fb^{\OI,d_i}_{\Pb}$.  By \cite[Theorem 6.15]{NR19}, $\Fb$ is a Noetherian OI-module. Hence, we conclude by using \Cref{cor:limit is noeth}. 
\end{proof}

In order to utilize systematically results on OI-modules for investigating Inc-modules we want to identify a given Inc-module as a colimit of an OI-module. This is indeed possible  for submodules of Noetherian free Inc-modules. As a substitute for the width of an element in an OI-module, we introduce the norm of an element of a free Inc-module that arises as a colimit. We use the notation of \Cref{rem:hom set as colimit} that associates to any $\pi \in \Hom_{\OI} (d, n)$ a map $\tilde{\pi} \in \Hom_{\OI} (d, \N)$ with the same image as $\pi$. 

\begin{defn}
    \label{def:norm}
Fix a finitely generated free Inc-module $\Fc = \bigoplus_{i = 1}^r \Fc^{\Inc, d_i}_{\Pc} $ with basis $\{\eps_{ \widetilde{\id_{d_i}}, i} \; \mid \; i \in [r]\}$ and integers $d_i \ge 0$ and consider an Inc-submodule $\Mcc$ with generating set $G = \{g_1,\ldots,g_s\}$, that is, $\Mcc = \langle G \rangle_{\Fc}$. 
Since $\Fc =  {\displaystyle \lim_{\longrightarrow}} \, \Fb$ with $\Fb = \bigoplus_{i = 1}^r \Fb^{\OI,d_i}_{\Pb} $ by \Cref{prop:colimit of free OI-module}, there are natural maps $\kappa_n \colon \Fb_n \to \Fc$. Define the \emph{norm}  of an element $q \in \Fc$ as the integer 
\[
\nu (q) = \min \{ n \in \N_0 \; \mid \; q \in \kappa_n (\Fb_n) \}. 
\]
The Inc-equivariant homomorphism 
\[
\phi \colon \Gc = \bigoplus_{i = 1}^s \Fc^{\Inc, \nu (g_i)}_{\Pc} \to \Mcc  \text{ defined by }  \gamma_{ \id_{\nu (g_i)}, i} \mapsto g_i, 
\]
where $\{\gamma_{\widetilde{ \id_{\nu (g_i)}}, i} \; \mid \; i \in [s]\}$ is a basis of $\Gc$, 
is called the \emph{surjection determined by $G$ and $\Fc$}. 
\end{defn}

If the free module $\Fc$ is clear from context, we refer to $\phi$ also simply as the surjection determined by $G$. 

\begin{ex}
   \label{exa:norm}
For a polynomial $q \in \Pc \cong \Fc^{\Inc, 0}_{\Pc}$, its norm is the maximum column index of a variable dividing a monomial in $q$. 
\end{ex}

Using the above notation, we are ready to identify $\Mcc$ as a colimit. 

\begin{prop}
    \label{prop:identify Inc-mod as limit}
If $K$ is a Noetherian ring then any submodule $\Mcc$ of a finitely generated free Inc-module over $\Pc$      
 is a colimit of an OI-module. In fact, one has $\Mcc = \lim \Mb$, where $\Mb$ is the submodule of $\Fb = \bigoplus_{i = 1}^r \Fb^{\OI, d_i}_{\Pb}$ that is generated by $b_i \in \Fb_{\nu(g_i)}$ with $\kappa_{\nu (g_i)} (b_i) = g_i$. 
\end{prop}

\begin{proof}
Note that by the definition of the norm $\nu (g_i)$ such an element $b_i \in \Fb_{\nu(g_i)}$ exists. It fact, it is unique (see \Cref{prop:colimit of free OI-module} and \Cref{exa:poly Inc-algebras}). Now one checks that $ \lim \Mb = \Mcc$. 
\end{proof}

\begin{rem}
If $\Ic$ is an ideal of $\Pc$, the construction of an OI-ideal $\Ib$ of $\Pb$ with $ \lim  \Ib = \Ic$. 
 can be made explicit and was known though not in this language (see \cite{HS, NR19}). In fact, $\Ib$ is the submodule of $\Pb$ defined by $\Ib_n = \Ic \cap K[X_{c \times n}]$, where $K[X_{c \times n}]$ is considered as a subring of $K[X_{c \times \N}]$. 

In general, there are several OI-ideals of $\Pb$ with the same colimit. Consider for example, the ideal $\Ic = \langle x_1^2, x_2^2, \ldots \rangle$ of $\Pc$ with $c =1$, where we write $x_j$ instead of $x_{1, j}$.  The above construction gives the OI-ideal $\Ib$ that is generated in width one by $x_1^2$. However, the OI-ideal $\Jb$ generated in width two by $x_1^2$ and $x_2^3$ has also $\Ic$ as colimit. 
\end{rem}

For any surjective morphism of Inc-modules over $\Pc$
\[
\phi \colon \Fc \to \Mcc, 
\]
where $\Fc$ is a free Inc-module, we call $\ker \, \phi$ a \emph{module of syzygies} of $\Mcc$. If $\Fc$ is finitely generated then \Cref{cor:fg free Inc-mod is noeth} shows that $\ker \phi$ is finitely generated. In fact, if $\Mcc$ is a submodule of finitely generated free Inc-module then a generating set of $\ker \phi$ can be determined algorithmically, provided $K$ is a field. 

\begin{proc}[Method for Computing Syzygies] Assume $K$ is a field.  
   \label{proc:Inc-syzygies}
   
Input: A finite subset $B$ of a finitely generated free Inc-module $\Fc = \bigoplus_{i = 1}^r \Fc^{\Inc, d_i}_{\Pc} $ with basis 
$\{\eps_{ \widetilde{\id_{d_i}}, i} \; \mid \; i \in [r]\}$. 

Output: An Inc-equivariant surjective homomorphism $\phi \colon \Gc 
\to \Mcc  = \langle B \rangle_{\Fc} $, where $\Gc$ is a finitely generated free Inc-module, 
and a finite generating set of $\ker \phi$. 
\begin{enumerate}

\item Find the surjection determined by $B$ and $\Fc$,  
\[
\phi \colon \Gc = \bigoplus_{i = 1}^s \Fc^{\Inc, \nu (b_i)}_{\Pc} \to \Mcc  \text{ defined by }  \gamma_{\widetilde{ \id_{\nu (b_i)}}, i} \mapsto b_i, 
\]
where $B = \{b_1,\ldots,b_s\}$ and $\{\gamma_{  \widetilde{\id_{\nu (b_i)}}, i} \; \mid \; i \in [s]\}$ is a basis of $\Gc$, 
as described in \Cref{def:norm}. 

\item Set $\Fb = \bigoplus_{i = 1}^r \Fb^{\OI, d_i}_{\Pb} $ with basis $\{e_{ \id_{d_i}, i} \; \mid \; i \in [r]\}$  
and $\Gb =   \bigoplus_{i = 1}^s \Fb^{\OI, \nu (b_i)}_{\Pb}$ with basis $\{g_{ \id_{\nu (b_i)}, i} \; \mid \; i \in [s]\}$. Suppose $B = \{b_1,\ldots,b_s\}$ with 
\[
b_i = \sum_{j = 1}^r a_{i, j} \eps_{\pi_{i, j}, j}. 
\]
Set 
\[ 
q_i = \sum_{j = 1}^r a_{i, j} e_{\rho_{i, j}, j} \in \Fb_{\nu (b_i)}, 
\]
where $\rho_{i, j} \in \Hom_{\OI} (d_j, \nu (b_i))$ is defined by $\rho_{i, j} (k)  = \pi_{i, j} (k)$. Note that, using the notation of \Cref{rem:hom set as colimit}, we have $\widetilde{\rho_{i, j}} = \pi_{i, j}$, and so $\Mcc = {\displaystyle \lim_{\longrightarrow}} \, \Mb$ with $\Mb = \langle q_1,\ldots,q_r \rangle_{\Fb}$ and $[q_i] = b_i$. 

\item Define a morphism of free OI-modules
\[
\Phi \colon \Gb \to \Fb  \quad \text{ by } \Phi (g_{ \id_{\nu (b_i)}, i} ) = q_i. 
\]
Use \Cref{lem:oikergen} to determine a finite generating set $\{p_1,\ldots,p_t \}$ of $\ker \Phi$. Since $\Gc = {\displaystyle \lim_{\longrightarrow}} \, \Gb$, it follows by 
 \Cref{prop:colomit of OI-mod} that $\{[p_1],\ldots,[p_t] \}$ is a generating set of $\ker \phi$. 

\end{enumerate}

\end{proc} 

\begin{rem}
If $\Fc$ is graded and $B$ consists of homogeneous elements then $\Mcc$ is graded and the grading of $\Gc$ can be adjusted  so that the surjection $\phi$ and the module $\ker \phi$ become graded.  
In fact, if $b_i \in B$ has degree $t_i$ then set $\Gc = \bigoplus_{i = 1}^s \Fc^{\Inc, \nu (b_i)}_{\Pc}(-t_i)$ and $\Gb = \bigoplus_{i = 1}^s \Fb^{\OI, \nu (b_i)}_{\Pb} (-t_i)$. The corresponding maps $\phi$ and $\Phi$ in Steps (i) and (iii) are then graded morphisms. 
\end{rem} 

Iterating the above procedure, one gets an algorithm to compute (truncated) free Inc-equivariant resolutions. The existence of such resolutions is true in greater generality.  

\begin{thm}
    \label{thm:Inc resolutions} 
If $K$ is a Noetherian ring then any submodule $\Mcc$ of a finitely generated free Inc-module over $\Pc = K[X_{c \times \N}]$      admits an $\Inc$-equivariant  free resolution 
\[
\cdots \to \Fc_2 \to  \Fc_1 \to \Mcc \to 0
\]
with $\Inc$-equivariant differentials and free $P$-modules $\Fc_i$ that are finitely generated as $\Inc$-modules. 

Moreover, if $\Mcc$ is a graded module then a graded such resolution exists, i.e., the differentials have degree zero. 
\end{thm}

\begin{proof}
By \Cref{prop:identify Inc-mod as limit}, $\Mcc$ is the colimit of a finitely generated OI-module $\Mb$ over $\Pb$. Hence, $\Mb$ admits a free resolution by finitely generated free OI-modules  due to \cite[Theorem 7.1]{NR19}. Passing to the colimit and using \Cref{prop:colomit of OI-mod}, we get the 
desired Inc-Invariant resolution. 
\end{proof}

If $K$ is a field such an Inc-equivariant resolution can be determined up to any desired finite homological degree by iterating \Cref{proc:Inc-syzygies}. However, it is computationally more efficient to rely mostly on computations for OI-modules as in the above proof.

\begin{proc}[Method for Computing Free $\Inc$-Resolutions] Assume $K$ is a field.  
   \label{proc:Inc-free res}
   
Input: A finite subset $B$ of a finitely generated free Inc-module $\Fc = \bigoplus_{i = 1}^r \Fc^{\Inc, d_i}_{\Pc} $ with basis 
$\{\eps_{ \widetilde{\id_{d_i}}, i} \; \mid \; i \in [r]\}$, $s \ge 1$ an integer.

Output: An exact sequence of finitely generated Inc-modules over $\Pc$, 
\[
 \Fc_s \to \Fc_{s-1} \to \cdots \to  \Fc_1 \to \Mcc \to 0
\]
with Inc-equivariant differentials,  $\Mcc = \langle B \rangle_{\Fc}$ and free Inc-modules $\Fc_i$. 
\begin{enumerate}

\item Use Steps (i) and (ii) of \Cref{proc:Inc-syzygies} to determine a finitely generated submodule $\Mb$ of a free OI-module $\Fb$ over $\Pb$ with $\Mcc = {\displaystyle \lim_{\longrightarrow}} \, \Mb$  and $\Fc = {\displaystyle \lim_{\longrightarrow}} \, \Fb$. 

\item Use \cite[Procedure 5.1]{MN} to determine the beginning of a free OI-resolution of $\Mb$ over $\Pb$, 
\[
 \Fb_s \to \Fb_{s-1} \to \cdots \to  \Fb_1 \to \Mb \to 0. 
\]

\item Pass to the colimit to get the desired Inc-resolution. 

\end{enumerate} 

\end{proc}

Using results on free resolutions of OI-ideals, we obtain resolutions of their limits. We illustrate this principle with a few examples. First, we consider ideals in $\Pc = K[X_{1 \times \N}] = K[X_\N]$, where there is only one row of variables. In this case, we use again simplified notation and denote the variables by $x_1,x_2,\ldots$.  

\begin{ex}
   \label{exa:Inc-Koszul}
For any integer $k \ge 1$, consider the Inc-ideal $\Ic \subset \Pc= K[X_\N]$ that is generated by $x_1^k$, that is, 
\[
\Ic = \langle x_i^k \; \mid \; i \in \N \rangle. 
\]
Thus, if $k =1$ then $\Ic$ is the unique homogeneous maximal ideal of $\Pc$. 

Now consider the OI-ideal $\Ib$ of $\Pb = X^{\OI, 1}$ that is generated by $x_1^k$ in width one. Its width $n$ component is 
\[
\Ib_n = \langle x_1^k,\ldots, x_n^k \rangle \subset \Pb_n. 
\]
It follows that $\Ic = \lim \Ib$. By \cite[Example 8.7]{NR19}, the ideal $\Ib$ has the following graded free OI-resolution
\[
\cdots \to \Fc^{\OI, j}_{\Pb}(-j k) \to \Fc^{\OI, j-1}_{\Pb}(-(j-1) k) \to \cdots \to \Fc^{\OI, 1}_{\Pb}(-k) \to \Ib \to 0. 
\]
Note that the width $n$ component of this resolution is the classical Koszul complex resolving $\Ib_n$ as graded $\Pb_n$-module. By \Cref{prop:colomit of OI-mod} and \Cref{prop:colimit of free OI-module}, passing to the colimit gives a graded Inc-invariant  resolution of $\Ic$ by finitely generated free Inc-modules: 
\[
\cdots \to \Fc^{\Inc, j}_{\Pc}(-j k) \to \Fc^{\Inc, j-1}_{\Pc}(-(j-1) k) \to  \cdots \to \Fc^{\Inc, 1}_{\Pc}(-k) \to \Ic \to 0. 
\]
\end{ex}

In the next example, we consider a squarefree monomial ideal. 

\begin{ex}
   \label{exa:Inc-square-free}
For any integer $k \ge 1$, let $\Jc \subset \Pc = K[X_\N]$ be the Inc-ideal that is generated by all squarefree monomials of degree $k$ in $\Pc$. Thus, 
\[
\Jc = \langle x_{i_1} x_{i_2} \cdots x_{i_k} \; \mid \; 1 \le i_1 < i_2 < \cdots < i_k \rangle. 
\]
We want to resolve $\Jc$ by finitely generated free Inc-modules. 

To this end, consider the OI-ideal $\Jb$ of $\Pb = \Pb^{\OI, 1}$ that is generated by $x_1 x_2 \cdots x_k$ in width $k$. For $n \ge k$, its width $n$ component is 
\[
\Jb_n = \langle x_{i_1} x_{i_2} \cdots x_{i_k} \; \mid \; 1 \le i_1 < i_2 < \cdots < i_k \le n \rangle  \subset \Pb_n. 
\]
Thus, we get  $\Jc = \lim \Jb$. 

By \cite[Corollary 4.11]{FN21},  the ideal $\Jb$ has a graded free OI-resolution of the form
\begin{align*}
\cdots \to( \Fc^{\OI, k+j-1}_{\Pb})^{\binom{k+j-2}{k-1}}(-k -j +1) \to & (\Fc^{\OI, k+j-2}_{\Pb})^{\binom{k+j-3}{k-1}}(-k- j+2) \to \cdots \\[2pt]
\to &  ( \Fc^{\OI, k+1}_{\Pb})^{k}(-k - 1) \to  \Fc^{\OI, k}_{\Pb}(-k) \to \Jb \to 0.
\end{align*}
Passing to the colimit, we get as above the desired Inc-invariant free resolution
\begin{align*}
\cdots \to( \Fc^{\Inc, k+j-1}_{\Pc})^{\binom{k+j-2}{k-1}}(-k -j +1) \to & (\Fc^{\Inc, k+j-2}_{\Pc})^{\binom{k+j-3}{k-1}}(-k- j+2) \to \cdots \\[2pt]
\to &  ( \Fc^{\Inc, k+1}_{\Pc})^{k}(-k - 1) \to  \Fc^{\Inc, k}_{\Pc}(-k) \to \Jc \to 0.
\end{align*}
\end{ex} 

Finally, we consider an ideal of a polynomial ring with any finite number of rows of variables. 

\begin{ex}
    \label{exa:Ferrers ideal}
For any integer $k \ge 1$, let $\Ic \subset  K[X_{k \times \N}] = \Pc$ be the  Inc-ideal that is generated by $x_{1,1} x_{2,2} \cdots x_{k,k}$. Thus, 
\[
\Ic = \langle x_{1, i_1} x_{2, i_2} \cdots x_{k, i_k} \; \mid \; 1 \le i_1 < i_2 < \cdots < i_k \rangle. 
\]    
We want to determine  an Inc-invariant graded free resolution of  $\Ic$ with finitely generated free Inc-modules. 

To this end, consider the OI-ideal $\Ib$ of $\Pb = (X^{\OI, 1})^{\otimes k}$ that is generated in width $k$ by $x_{1,1} x_{2,2} \cdots x_{k,k}$. For $n \ge k$, its width $n$ component is 
\[
\Ib_n = \langle x_{1, i_1} x_{2, i_2} \cdots x_{k, i_k} \; \mid \; 1 \le i_1 < i_2 < \cdots < i_k \le n \rangle, 
\]   
and so $\lim \Ib = \Ic$. 

Note that, by \cite[Theorem 4.10]{FN21},  the free OI-resolution of the ideal $\Jb$ in \Cref{exa:Inc-square-free} is a cellular resolution supported on a cell complex. According to \cite[Theorem 4.13]{FN21}, the same cell complex also supports a graded free OI-resolution of the ideal $\Ib$. Thus, the ideal $\Ib$ has a free OI-resolution of the form 
\begin{align*}
\cdots \to( \Fc^{\OI, k+j-1}_{\Pb})^{\binom{k+j-2}{k-1}}(-k -j +1) \to & (\Fc^{\OI, k+j-2}_{\Pb})^{\binom{k+j-3}{k-1}}(-k- j+2) \to \cdots \\[2pt]
\to &  ( \Fc^{\OI, k+1}_{\Pb})^{k}(-k - 1) \to  \Fc^{\OI, k}_{\Pb}(-k) \to \Ib \to 0.
\end{align*}
Passing to the colimit, we get as above the desired Inc-invariant free resolution
\begin{align*}
\cdots \to( \Fc^{\Inc, k+j-1}_{\Pc})^{\binom{k+j-2}{k-1}}(-k -j +1) \to & (\Fc^{\Inc, k+j-2}_{\Pc})^{\binom{k+j-3}{k-1}}(-k- j+2) \to \cdots \\[2pt]
\to &  ( \Fc^{\Inc, k+1}_{\Pc})^{k}(-k - 1) \to  \Fc^{\Inc, k}_{\Pc}(-k) \to \Ic \to 0.
\end{align*}
\end{ex}


\section{Colimits of FI-modules}
\label{section:Sym-resolutions} 

We introduce $\Sym$-modules over $\Sym$-algebras, where 
\[
\Sym = \{\pi \colon \N \to \N \; \mid \; \pi \text{ is bijective and $\pi (i) \neq i$ for finitely many } i \in \N \}. 
\]
It turns out that colimits of FI-modules are $\Sym$-modules. Using our results about FI-resolutions, we establish an algorithm for computing  a $\Sym$-equivariant free resolution of any finitely generated $\Sym$-module over $K[X_{c \times \N}]$, considered as a $\Sym$-algebra, up to any finite homological degree. In particular, it applies to any symmetric ideal of $K[X_{c \times \N}]$. 
 Most of the results and proofs are similar to statements in the previous section, and so we mainly restrict ourselves to mentioning notable adjustments or differences. 

Note that $\Sym$ can be identified with $\bigcup_{n \in \N} \Sym (n)$, where $\Sym(n)$ is embedded into $\Sym(n+1)$ as the stabilizer of $n+1$. Though we will always consider elements of $\Sym$ as maps $\N \to \N$. 

Denote by $K$ any commutative ring. Analogously to \Cref{def:Inc-algebra}, we define a \emph{$\Sym$-algebra} over $K$ as a  $K$-algebra $\Ac$ that admits an action $\Sym \times \Ac \to \Ac, \tau \times a \mapsto \tau \cdot a$,  such that each induced map $\tau(-) \colon \Ac \to \Ac$ is a $K$-algebra homomorphism.  A \emph{$\Sym$-algebra homomorphism} is a $\Sym$-equivariant homomorphism of $\Sym$-algebras over $K$. Thus, $K[X_{c \times \N}]$ is a $\Sym$-algebra with an action induced by $\tau \cdot x_{i, j} = x_{i, \tau (j)}$. 

For any integer $d \ge 0$,  there  is a natural action of  $\Sym$  on  the set 
\[
\Hom_{\FI} (d, \N) = \{ \pi \colon d \to \N \; \mid \; \text{ $\pi$  is injective} \}, 
\]  
defined by $\tau \times \pi = \tau \circ \pi$. Given any $\tau \in \Sym$ and $\pi \in \Hom_\FI (d, n)$, set $k = \max \{ \tau(\pi (j)) \; \mid \; j \in [d]\}$ and define a map 
\[
\tau \cdot \pi \colon [d] \to [k], \; j \mapsto \tau (\pi (j)). 
\] 
It is in $\Hom_\FI (d, k)$. Then, for any $\FI$-algebra $\Ab$ over $K$, using the direct system $(\Ab_n, \Ab(\iota_{m, n}))$ with maps $\iota_{m, n}$ introduced in \Cref{rem:hom set as colimit},  its colimit $\Ac = {\displaystyle \lim_{\longrightarrow}} \, \Ab_n$ is a $\Sym$-algebra whose $\Sym$-action is defined by
\[
\tau \cdot [a] = [\Ab(\tau \cdot \id_n) (a)], 
\]
where $[a]$ denotes the class of $a \in \Ab_n$ and $\tau \in \Sym$. For example, for $\Pb^{\FI, c}$ we get as colimit 
\[
{\displaystyle \lim_{\longrightarrow}} \, \Pb^{\FI, c} \cong K[X_{c \times \N}]  
\]
where $K[X_{c \times \N}]$ is considered as $\Sym$-algebra with the above $\Sym$-action and the class $[f]$ of any $f \in \Pb^{\FI, c}_n$ is mapped onto $f$, considered as a polynomial in $K[X_{c \times \N}]$. 

For any $\Sym$-algebra $\Ac$ over $K$, define the category of \emph{$\Sym$-modules over $\Ac$}, denoted $\Sym$-$\Mod (\Ac)$, as in \Cref{def:Inc-module} by replacing $\Inc$ by $\Sym$. It is an abelian category. 
Similarly, a \emph{free} $Sym$-module over $\Ac$ is any $\Sym$-module $\Fc$ that is isomorphic to  direct sum $\bigoplus_{\lambda\in\Lambda}\Fc^{\Sym,d_\lambda}_{\Ac} $ with integers $d_\lambda\geq0$, where, for any $d \ge 0$,  
\[
\Fc^{\Sym,d} = \Fc^{\Sym,d}_{\Ac} =\bigoplus_{\pi\in\Hom_{\FI} (d,\N)}\Ac \eps_\pi   
\]
is a free $\Ac$-module with basis $\{ \eps_{\pi} \; \mid \; \pi \in \Hom_{\FI} (d,\N) \}$ 
and action $\Sym \times \Fc^{\Inc,d} \to \Fc^{\Inc,d}$ induced by
\[
\tau \cdot (a e_{\pi}) = (\tau \cdot a) e_{\tau \circ \pi}. 
\]
If $|\Lambda|=r<\infty$, then $\Fc$ is said to have \emph{rank $r$}. Observe that  $\Fc^{\Sym,d}$ is generated by one element as a 
$\Sym$-module, namely by  $\eps_{\widetilde{\id_d}}$, where $\widetilde{\id_d}\colon [d] \to \N, j \mapsto j$. 

Analogously to \Cref{prop:colomit of OI-mod}, one shows that, for any $\FI$-algebra $\Ab$,  passing to colimits gives an exact functor from the category of $\FI$-modules over $\Ab$, to the category of $\Sym$-modules over $\Ac = {\displaystyle \lim_{\longrightarrow}} \, \Ab$. In particular, for any $d \ge 0$,  there is an isomorphism of $\Sym$-modules 
\[
\lim_{\longrightarrow} \Fb^{\FI,d}_{\Ab} = \Fc^{\Sym,d}_{\Ac}. 
\]
If $K$ is a Noetherian ring, then it follows as in the previous section that any finitely generated $\Sym$-module over the $\Sym$-algebra $K[X_{c \times \N}]$ is Noetherian. 
 
For any finitely generated free $\FI$-module $\Fb$, define the \emph{norm} of an element 
$q \in \Fc =  {\displaystyle \lim_{\longrightarrow}} \, \Fb$ as in \Cref{def:norm}.  With these preparations in place, one obtains the following analog of \Cref{thm:Inc resolutions}. 

\begin{thm}
    \label{thm:Sym resolutions} 
If $K$ is a Noetherian ring then any submodule $\Mcc$ of a finitely generated free $\Sym$-module over $\Pc = K[X_{c \times \N}]$      admits a $\Sym$-equivariant  free resolution 
\[
\cdots \to \Fc_2 \to   \Fc_1 \to \Mcc \to 0
\]
with $\Sym$-equivariant differentials and free $P$-modules $\Fc_i$ that are finitely generated as $\Sym$-modules. 

Moreover, if $\Mcc$ is a graded module then a graded such resolution exists, i.e., the differentials have degree zero. 
\end{thm} 

If $K$ is a field, such a $Sym$-equivariant free resolution can be determined up to any desired finite homological degree in finitely many steps by using \Cref{proc:freefires} and an analog of \Cref{proc:Inc-free res}. 

\begin{proc}[Method for Computing Free $\Sym$-Resolutions] Assume $K$ is a field.  
   \label{proc:Sym-free res}
   
Input: A finite subset $B$ of a finitely generated free $\Sym$-module $\Fc = \bigoplus_{i = 1}^r \Fc^{\Sym, d_i}_{\Pc} $ with basis 
$\{\eps_{ \widetilde{\id_{d_i}}, i} \; \mid \; i \in [r]\}$, $s \ge 1$ an integer.

Output: An exact sequence of finitely generated $\Sym$-modules over $\Pc$, 
\[
 \Fc_s \to \Fc_{s-1} \to \cdots \to  \Fc_1 \to \Mcc \to 0
\]
with Inc-equivariant differentials,  $\Mcc = \langle B \rangle_{\Fc}$ and free Inc-modules $\Fc_i$. 
\begin{enumerate}

\item Similarly to Steps (i) and (ii) of \Cref{proc:Inc-syzygies}, determine a finitely generated submodule $\Mb$ of a free FI-module $\Fb$ over $\Pb$ with $\Mcc = {\displaystyle \lim_{\longrightarrow}} \, \Mb$  and $\Fc = {\displaystyle \lim_{\longrightarrow}} \, \Fb$. 

\item Use \Cref{proc:freefires} to determine the beginning of a free FI-resolution of $\Mb$ over $\Pb$, 
\[
 \Fb_s \to \Fb_{s-1} \to \cdots \to  \Fb_1 \to \Mb \to 0. 
\]

\item Pass to the colimit to get the desired Sym-resolution. 

\end{enumerate} 

\end{proc}


\end{document}